\documentclass[11pt]{amsart}
\usepackage{amsmath,amsthm,amsfonts,amssymb}

\usepackage{amsfonts}\usepackage{amsthm}
\usepackage{enumerate}
\usepackage{graphicx}

\def\leaderfill{\leaders\hbox to 1em{\hss.\hss}\hfill}

\theoremstyle{plain}
\newtheorem{thm}{Theorem}[section]

\newtheorem{lem}[thm]{Lemma}
\newtheorem{lemma}[thm]{Lemma}

\newtheorem{prop}[thm]{Proposition}
\numberwithin{equation}{section}

\theoremstyle{definition}

\theoremstyle{remark}

\numberwithin{equation}{section}

\def\volume{\text{area}}
\def\a{{\alpha}}

\def\vth{\vartheta}

\def\cf{{\it cf}}

\def\diam{{\text{\rm diam}}}
\def\area{{\text{\rm area}}}

\def\supp{{\text{\rm supp }}}

\def\inn#1#2{\langle#1,#2\rangle}

\def\noi{\noindent}

\def\meas{{\text{\rm meas}}}

\def\card{\text{\rm card}}
\def\lc{\lesssim}
\def\gc{\gtrsim}
\def\ga{\gamma}             
\def\eps{\varepsilon}

\def\ep{\epsilon}

\def\la{\lambda}

\def\om{\omega}              \def\Om{\Omega}
\def\ka{\kappa}

\def\bbE{\mathbb{E}}
\def\bbG{\mathbb{G}}
\def\bbA{\mathbb{A}}
\def\bbB{\mathbb{B}}
\def\bbI{\mathbb{I}}
\def\bbJ{\mathbb{J}}
\def\bbR{\mathbb{R}}
\def\bbZ{\mathbb{Z}}

\def\bbT{\mathbb{T}}
\def\cA{\mathcal{A}}

\def\cC{\mathcal{C}}

\def\cG{\mathcal{G}}

\def\cI{\mathcal{I}}

\def\cK{\mathcal{K}}
\def\cL{\mathcal{L}}

\def\cN{\mathcal{N}}

\def\cS{\mathcal{S}}
\def\cT{\mathcal{T}}

\def\ic{i}

\def\fS{{\mathfrak S}}

\def\fB{{\mathfrak B}}

\def\fP{{\mathfrak P}}
\def\fT{{\mathfrak T}}
\def\fL{{\mathfrak L}}
\def\cD{\mathcal{D}}

\def\fn{{\mathfrak n}}
\def\kunit{{\tfrac{k}{|k|}}}
\def\lunit{{\tfrac{\ell}{|\ell|}}}

\begin{document}

\title[Mean lattice point discrepancy bounds]
{
 Mean lattice point discrepancy bounds, II:
Convex domains in the plane}
\author{Alexander Iosevich \ \  Eric T. Sawyer  \ \  Andreas Seeger}
\address{
A. Iosevich \\
Mathematics Department\\ University of Missouri\\
 Columbia, MO 65211, USA}
\email{ iosevich@math.missouri.edu}
\address{ E.T. Sawyer\\
Department of Mathematics and Statistics\\
 1280 Main Street West\\
 Hamilton, Ontario  L8S 4K1, Canada}
\email{sawyer@mcmaster.ca}
\address{
 A. Seeger\\
 Department of
 Mathematics\\
University of Wisconsin-Madison\\
Madison, WI 53706, USA}
\email{seeger@math.wisc.edu}
\thanks{Research supported in part by  NSF grants.}
\date {Revised version, January 15, 2007
(originally submitted version: September 26, 2004).}

\begin{abstract}
We consider planar curved strictly convex domains with no or very weak
smoothness assumptions and prove sharp bounds  for
square-functions associated to the lattice point discrepancy.
\end{abstract}

\maketitle


\section{Introduction}
 This  paper is a sequel to \cite{ISS} in which the authors proved
bounds for the mean square lattice point discrepancy for convex bodies
with smooth boundary
in ${\mathbb R}^d$. Here we reconsider the case $d=2$
but admit now domains with rough boundary.

Let $\Omega$ be a convex domain in $\bbR^2$ containing
the origin in its interior.
Let $$\cN_\Om(t)=\card(t\Omega\cap\bbZ^2),
$$
the number of integer lattice points inside the dilate $t\Omega$.
It is well known  that $\cN_\Om(t)$ is asymptotic
to $t^2\text{area}(\Om)$ as $t\to \infty$ and we denote by
\begin{equation}\label{eq:1.1}
E_\Om(t)=\cN_\Om(t)-t^2\volume(\Om)
\end{equation}
the  error, or  lattice rest.
A trivial estimate for the lattice rest is
$E_\Om(t)\le Ct$ which holds for any convex set.
For the case that the boundary is smooth and has positive curvature everywhere
this estimate has been significantly improved. It is  conjectured that in this case
$E_\Om(t)=O(t^{1/2+\eps})$ for any $\eps>0$ but by the best
result  published at this time,  due
to Huxley \cite{Hu2},
one only knows that
$E_\Om(t)= O( t^{131/208}(\log t)^{A})$
for suitable $A$.

On average however better estimates hold.
We consider the mean-square discrepancy of the lattice rest over the interval
$[R, R+h]$  where $h\le R$ and $R$ is large; it is given by
\begin{equation}\label{eq:1.2}
\cG_\Om(R,h)
=\Big(\frac 1h \int_R^{R+h}|E_\Om(t)|^2 dt\Big)^{1/2}.
\end{equation}

Provided that the boundary is smooth (say $C^4$) and the Gaussian curvature
never vanishes
it has been shown by Nowak \cite{No1} that $\cG_\Om(R,R)=O(R^{1/2})$;
later
Huxley \cite{Hu} showed that
$\cG_\Om(R,1)=O(R^{1/2}\log^{1/2}R)$.
A result which unifies both estimates is in the authors' paper
\cite{ISS}, namely
$\cG_\Om(R,h)\le CR^{1/2}$ if $\log R\le h\le R$. We
 note that Nowak \cite{No2} has independently proved the same bound.
Moreover he
obtained asymptotics for $\cG_\Om (R, H(R))$ as $R\to \infty$, provided that
$H(R)/\log R\to \infty$; see also  earlier asymptotics by
Bleher \cite{B} where essentially $H(R)\approx R$.

The purpose of this paper is
 to prove versions of these estimates under minimal
(or no) smoothness assumptions
on the boundary of the domain.  The main
 difficulty is that the oscillation of the Fourier transforms of densities on
the boundary
cannot be used in a straightforward way as in \cite{No}, \cite{ISS},
 or \cite{No2},
 because of the lack of  asymptotic expansions.

Our first result deals with
 domains for which the curvature is bounded below
with very  weak regularity assumptions on the curvature.
Here we assume that $\Omega$ has $C^1$ boundary,
that the components of the tangent vectors are absolutely continuous functions
of the arclength parameter so  that the second derivatives of a regular parametrization
are well defined as   $L^1$ functions on the boundary. The following theorem
yields an analogue of the above result with  a slightly
more restrictive assumption on these second derivatives.

\begin{thm}\label{thm:1.1} Let $\Omega$
be a convex domain in $\bbR^2$ containing the
origin in its interior, and assume that  $\Omega$ has  $C^{1}$  boundary and that the
components of the tangent vector are absolutely continuous functions.
 Suppose also that
 curvature $\kappa$ is bounded below, i.e. $\kappa(x)\ge a>0$ for almost every $x\in\partial \Om$
and that $\kappa\in L\log^{2+\epsilon}L(\partial\Omega)$, for some $\epsilon>0$.
Then there is a constant $C_\Om$ so that for all $R\ge 2$
\begin{equation}
\label{eq:1.3}
\cG_\Om(R,h)\le
 C_\Om
R^{1/2} \quad \text{ if }\quad \log R\le h\le R.
\end{equation}
\end{thm}

Of course this result applies to all convex domains with $C^2$ boundary
and nonvanishing curvature; but it also applies to rougher domains, the
simplest examples are $\{x:|x_1|^{a_1}+|x_2|^{a_2}\le 1\}$ when
$1<a_1,a_2\le 2$.  Moreover if ${\mathcal D}$ is a convex domain with smooth
finite type boundary, containing the origin, then the polar set
$\Omega=\cD^*=\{x:\sup_{\xi\in {\mathcal D}}\inn{x}{\xi}\le 1\}$ satisfies
the
 assumptions of Theorem \ref{thm:1.1}. For these examples the second
derivatives belong to
$L^p(\partial \Omega)$ for some $p>1$
 ({\it cf.}  the calculations in the proof of Lemma 5.1
 in \cite{ISS}).

An immediate consequence of (\ref{eq:1.3}) is Huxley's bound (\cite{Hu}) who proved that
$\cG_\Omega(R,1)=
O(\sqrt{R\log R})$ under the assumption that $\Om $ has $C^4$ boundary and the
curvature is bounded below.
We shall see ({\it cf.}  Theorem \ref{thm:1.3} below) that it  is possible
to prove this estimate for
convex domains  in which even  the weak  regularity assumption of Theorem
\ref{thm:1.1} is removed. Moreover in this case we shall prove ({\it cf.}
Theorem
 \ref{thm:1.2} below) that the optimal bound $\cG_\Omega(R,h)=O(R^{1/2})$
holds
in the more restricted range of $h$'s
 $R^{1/2}\le h\le R$.

In this rough case the assumption of the curvature bounded below has
to be reformulated
(as now we are not actually assuming that the curvature is a well defined  function).
Let $\rho^*$ be the Minkowski functional of the polar set $\Omega^*$, {\it i.e.}
\begin{equation}\label{eq:1.4}
\rho^*(\xi)=\sup\{\inn x\xi: x\in \Omega\}
\end{equation}
so that $\Omega^*=\{ \xi:\rho^*(\xi)\le 1\}$. For
$\theta\in S^1$ and $\delta>0$ consider the arc (or ``cap'')
\begin{equation}\label{eq:1.5}
\cC(\theta,\delta)\equiv \cC_\Omega(\theta,\delta)=\{x\in \partial
\Omega:\inn{x}{\theta}=\rho^*(\theta)-\delta\}.
\end{equation}
Let
\begin{equation}\label{eq:1.6}
\mu(\theta,\delta)=
\text{{\rm diam}}(\cC(\theta,\delta)).
\end{equation}
We note that if $d\sigma$ is the arclength measure on $\partial \Om$ then
  $\sum_\pm\mu(\pm\theta,\delta)$ controls the size of the Fourier
transform $\widehat {d\sigma }(\pm \theta/\delta)$, see \cite{BNW} and
also
\cite{BRT}. If the curvature is absolutely continuous and bounded below
then
it is easy to see that $\mu(\theta,\delta)=O(\sqrt\delta)$ uniformly in
$\theta\in S^1$, and in the general case we shall simply assume the
validity of this inequality.
It is possible to show the
equivalence of this condition with
 with other natural definitions of bounded below curvature for rough
 domains but we shall not discuss this here.

\begin{thm}\label{thm:1.2}
Let $\Omega$
be a convex domain in $\bbR^2$ containing the
origin in its interior. Suppose that
\begin{equation}\label{eq:1.7}
\sup_{\theta\in
S^1}\sup_{\delta>0}\delta^{-1/2}\mu(\theta,\delta)<\infty.
\end{equation}
Then for $R\ge 2$,
\begin{equation}\label{eq:1.8}
\cG_\Om(R,h)\le C_{\Om}
R^{1/2} \quad \text { if }\quad R^{1/2} \le h\le R.
\end{equation}
\end{thm}

If we admit an additional factor of $\sqrt{\log R}$ the range of $h$ can
be vastly improved to obtain a version of Huxley's theorem (\cite{Hu})
for rough domains with nonzero curvature (which
is much more elementary than
Theorem 1.2).

\begin{thm}\label{thm:1.3}
Let $\Omega$ be as in Theorem \ref{thm:1.2} (satisfying (\ref{eq:1.7})).
Then for $R\ge 2$
\begin{equation}\label{eq:1.9}
\cG_\Om(R,h)\le C_\Om
(R \log R)^{1/2} \quad \text { if }\quad 1\le h\le R.
\end{equation}
\end{thm}

\noindent {\it Remark}.
 An examination of the proof of  Theorem \ref{thm:1.3} shows that
the
constants depend only on the bound in (\ref{eq:1.7}) and
the radii of inscribed and circumscribed circles centered  at the origin.
This uniform version of inequality
(\ref{eq:1.9}) as well as the statement of Theorem  \ref{thm:1.1}
 is close to sharp as one can show that they  fail for $h\le
(\log R)^{-1}$.
To see this one uses Jarn\'\i k's   curve (\cite{J}) to produce a
sequence $R_j\to \infty $
and domains
 $\Omega_j$, so that the maximal inscribed and minimal  circumscribed radii
of $\Omega_j$ are  bounded above and
below,  the curvature on the boundary is bounded above and below and
$E_{R_j \Om_j}\ge R_j^{2/3}$ (\cite{J}, \cite{L}).
By  Huxley's mean-max inequality (\cite{Hu1}, p.\ 136)
$$\Big(\frac 1{\delta}
\int_{R_j-\delta}^{R_j+\delta}E_{\Om_j}(s)^2 ds\Big)^{1/2}\ge
E_{\Om_j}(R_j)/2$$
which holds under the assumptions that $|E_{\Om_j}(R_j) |
\ge 5(\text{area}(\Om_j))
\delta R_j$
and $0<\delta\le R_j/2$. We apply this for $\delta\approx R_j^{-1/3}\le h$
to see that under the assumption of
$\cG_{\Om_j}(R_j,h)
\lc (R_j\log R_j)^{1/2}$ we have
\footnote{Here, and in what follows
we use the following notation:
Given two quantities $A$, $B$ we write $A\lesssim B$
if there is an
absolute positive constant,
depending only  on the specific domain $\Omega$, so that $A\le C B$.
We write $A\approx B$ if $A\lc B$ and $B\lc A$.
}
$$R_j^{2/3}\lc E_{\Om_j}(R_j)\lc R_j^{1/6} h^{1/2} \cG_{\Om_j}(R_j,h)
\lc R_j^{2/3}(h\log R_j)^{1/2}.
$$
thus $h\gc (\log R_j)^{-1}$.
{\it Cf.} also Plagne \cite{Pl} for the construction of a single strictly convex curve $C$
  and a sequence $R_j$ so that $R_j C$ contains $R_j^{2/3} w(R_j)$ lattice
points, with $w(R)$ converging to zero at a slow rate.

It is certainly conceivable that the result of Theorem 1.2 may
hold for some $h\ll R^{1/2}$.  However this could not be established by simple extensions of our method, see the discussion below and in
\S \ref{sharpness}.

\medskip

Finally, if we consider arbitrary convex domains (dropping
 the curvature assumptions on the boundary) then the
estimate (\ref{eq:1.9}) may  fail
as does the classical estimate
$E_{\Om}(t)=O(t^{2/3})$ ({\it cf.} \cite{R}).
However for almost all rotations $A\in SO(2)$ it is still true that
$E_{A\Om}(t)=O(t^{2/3}\log^{1/2+\eps}t)$, see \cite{BCI}.
In fact for domains with smooth
finite type
boundary  one has the better estimate $E_{A\Om}(t)=O(t^{2/3-\delta})$,  for almost all rotations,
for some $\delta>0$,
 see Nowak's article \cite{No1}. Likewise  for such domains it is
proved in \cite{ISS} that for almost all rotations $A$
we have $\cG_{A\Om}(R,h)\lc R^{1/2}$ for all $R\ge 2$, $\log R\le h\le R$.
For arbitrary convex domains we can prove an analogous result but lose
 an additional  power of a logarithm.

\begin{thm}\label{thm:1.4}
Let $\Omega$
be a convex domain in $\bbR^2$ containing the
origin in its interior.
For $\vartheta \in [-\pi,\pi]$ denote by $A_\vartheta$ the rotation by the angle $\vartheta$ and by  $A_\vartheta \Omega$ the rotated
domain $ \{A_\vartheta x:x\in\Omega\}$.
Then for $\epsilon>0$, $R\ge 2$,
\begin{equation}\label{eq:1.10}
\cG_{A_\vartheta \Om}(R,h)\le
\cC_{\eps,\Om}(\vartheta)
R^{1/2} (\log R)^{1+\epsilon} \quad \text { if }\quad 1\le h\le R
\end{equation}
where $\cC_{\epsilon,\Om}(\vartheta )<\infty$ for almost all
$\vartheta\in [-\pi,\pi]$; in fact the function
$\cC_{\epsilon,\Om}$ belongs to the weak type space $ L^{2,\infty}$.
\end{thm}

{\it Structure of the paper:}
 The first parts of the proofs are  identical for
Theorems \ref{thm:1.1}--\ref{thm:1.4}. One uses essentially a
``$T^*T$-argument'' to reduce to a weighted
estimate for lattice points in thin annuli formed by dilations of the
polar domain.
This argument is straightforward in the smooth case (\cite{ISS},
\cite{No2})
but there are considerable  technical complications in the rough case.
The relevant estimates are given in \S2. One is
led to the estimation of quantities such as
\begin{equation}\label{bilinear}
\cK(R,h)=
\sum_{(k,\ell)\in \bbZ^2\times \bbZ^2
\atop {
{|k|,|\ell|\le R}
\atop {|\rho^*(k)-\rho^*(\ell)|\le h^{-1}}}
}
|k|^{-2} \mu(\kunit, \tfrac {1}{|k|R})\,
\mu(\lunit, \tfrac {1}{|\ell|R})
\end{equation}
when $h\gg1$,
and some variants with tails. \S3 contains further
discussion of these quantities in the case of curvature
bounded below and the main technical propositions needed
for the proofs of Theorems
\ref{thm:1.1}, \ref{thm:1.2}, and \ref{thm:1.3}.

The estimation of \eqref{bilinear} which is needed for
Theorems \ref{thm:1.3} and \ref{thm:1.4} is rather straightforward.
For Theorem 1.3 one uses the bound
$\mu(\theta,\delta)= O(\delta^{1/2})$
in conjunction with  the trivial bound
$E_{\Om^*}(t)=O(t)$.
For Theorem 1.4 one argues similarly but uses an averaged estimate  for
the size of caps for  the rotated domains.
The proofs are contained in \S3 and \S4.
The mild
regularity  assumption in Theorem \ref{thm:1.1} can be used to
improve on the trivial bound for $E_\Om^*(t)$. This
leads to boundedness
of $\cK(R,h)$ for $h\ge \log R$
and then in this range to the optimal  bound (\ref{eq:1.3});
the argument  is carried out in \S5.

The main technical estimate needed for the proof
of Theorem \ref{thm:1.2} is stated as Proposition \ref{prop:3.3}.
Here
 we need to efficiently estimate a {\it weighted}  version of the lattice
point discrepancy for $\Omega^*$ (\cf. \eqref{weightedsum}), and we shall
use  a more geometrical approach  for which we need the assumption
$h\ge R^{1/2}$. The proof is carried out in
\S6-\S8.

We do not know whether
in the generality of Theorem \ref{thm:1.2}
the assumption $h\ge R^{1/2}$ is really necessary.
In \S9 we construct some examples
of sets with rough boundary (and curvature bounded below)
which show that at least for the estimation of $\cK(R,h)$
the condition $h\ge R^{1/2}$ is necessary
 (which only shows the sharpness of the method).

Finally, a short appendix (\S\ref{appendix2}) is included which is not closely
 connected to the rest of the paper but  is relevant
to our previous results in \cite{ISS}. We note that following  
the  submission of the
current paper the arguments in this appendix
have already been used again in  subsequent work (\cite{IR}).
Here we discuss
a connection between mean discrepancy results and
 generalized distance sets
for  integer point lattices and by our previous results  obtain
 new lower bounds  in  three and higher dimensions.


\section{The first step}

Our purpose here is to show an estimate involving the
quantities $\mu(\theta, \delta)$, which holds without
any regularity or curvature assumptions on the boundary of the convex domain.

However we shall first make the {\it a priori}

\medskip

\noi{\bf Assumption:} The boundary  of $\Omega$ is a $C^1$
curve and the
 components of the outer  unit normal vectors  are absolutely continous.

\medskip

This means we assume that the curvature is integrable. Below we shall remove this {\it a priori} assumption by a limiting argument.

We now begin with a standard procedure using mollifiers to
regularize the characteristic function of $\Omega$. Suppose that
 $r_1< 1<r_2$
and $r_1,r_2$
 are the radii of inscribed and circumscribed circles centered at the origin.
Let $\zeta$ be a smooth nonnegative radial
 cutoff function supported in the ball $B_{r_1/2}(0)$
 so that $\int\zeta(x) dx=1$ and let $\zeta_{\eps}(x)=\eps^{-2}\zeta(x/\eps)$.
Let
\begin{equation}\label{Neps}
N_\eps(t)=\sum_{k\in \bbZ^2} \chi_{t\Omega} * \zeta_{\eps} (k)
\end{equation}
and
\begin{align}\label{Eeps}
E_\eps(t)&=\sum_{k\in \bbZ^2}
\chi_{t\Omega} * \zeta_{\eps} (k)- t^2 \volume(\Omega).
\end{align}

It suffices to estimate the modified square-function
\begin{equation}\label{squarefctmodified}
G(R,h)=\Big(\frac 1{h}\int_{R}^{R+h}|E_{1/R}(t)|^2 dt\Big)^{1/2 }
\end{equation}
where $h\ge 1$ 
since there is the elementary estimate (with $E\equiv E_\Omega$)
\begin{equation}\label{Emollified}
\Big(\frac 1{h}\int_{R}^{R+h}|E(t)|^2 dt\Big)^{1/2}\le
G(R,h)+C (R/h)^{1/2},
\end{equation}
valid for
$R^{-1}\le h\le R$;
see Lemma (2.2) of \cite{ISS}.

\medskip

\noi{\bf Basic decompositions.}
Fix a nonnegative  $ \eta_0\in C^\infty(\mathbb{R})$   so that  $\eta_0(t)
=1$ for
$t\in [0,1]$ and $\eta_0$ is supported in $ (-1/2,3/2)$ and let
\begin{equation}\label{etaRk}
\eta_{R, h}(t)=\frac{1}{\sqrt
h}\eta_0(\frac{t-R}{h}) .
\end{equation}
Then
\begin{equation}\label{etaRk2}
\frac 1h\int_R^{R+h}|E_{1/R}(t)|^2 dt\lc
\int |E_{1/R}(t)\eta_{R,h}(t)|^2 dt.
\end{equation}

By the Poisson summation formula
\begin{eqnarray}\label{poisson}
E_{1/R} (t)&=&\sum_{k\neq 0} (2\pi t)^2\widehat{\chi_\Om}
(2\pi tk)\widehat\zeta (2\pi k/R)
\\
&=&\sum_{0<|k|\le R^2} (2\pi t)^2\widehat{\chi_\Om}
(2\pi tk)\widehat\zeta (2\pi k/R) + O(R^{-10})\nonumber
\end{eqnarray}
since always
$|\widehat{\chi_\Om}
(2\pi tk)|\lc |tk|^{-1}$ and $|\widehat\zeta (2\pi k/R)|\le C_N(1+|k|/R)^{-N}$.

As in \cite{Hl} and elsewhere
we have by the divergence theorem
$\widehat {\chi_\Omega}(\xi)=-\ic\sum_{i=1}^2
(\xi_i/|\xi|^2)\widehat{\fn_i d\sigma}(\xi)$ where $\fn$ denotes the unit
outer normal vector. We
may  assume that the boundary of
$\Omega$ is parametrized by
$\alpha\mapsto x(\alpha)$ where $x'(\alpha)$ is a unit vector and
$x(\alpha)=x(\alpha+L)$ if $L$ is the length of $\partial \Om$.

 Then
$\fn(\a)=-x_\perp'(\a)$ where $x_\perp(\a)=(x_2(\a),-x_1(\a))$
and
\begin{equation}\label{div}
\widehat {\chi_\Omega}(\xi)= -\ic
\int_{0}^{L}
\frac{\inn{\xi}{x_\perp'(\a)}}{|\xi|^2} e^{-\ic \inn{x(\a)}{\xi}} d\a.
\end{equation}

Assuming that $R/2\le t\le 2R$ we shall now introduce a finer
microlocal
decomposition of $\widehat{\chi_\Om}(2\pi tk)$, depending on $k$ and $R$ and
 based on (\ref{div}).
This is  somewhat
inspired by \cite{BNW}, \cite{NSW} and in  particular by
 \cite{SZ} where a related construction  is used.
\medskip

 Suppose that $\beta_0$ is an even function which is supported in
$(-3/4,3/4)$ and which is equal to one in $[-1/2,1/2]$. Let
$\beta(s)=\beta_0(s/2)-\beta_0(s)$ and let, for $n\ge 1$,
$\beta_n(s)=\beta(2^{-n}(s))$.
Let

\begin{multline*}
\Psi(k,\alpha)=\frac{\inn{k}{x_\perp'(\a)}}{|k|^2}\quad\times
\\
\Big(1-\beta_0\big(2r_1^{-1}(
\inn{\tfrac{k}{|k|}}{x(\a)}-\rho^*(\tfrac{k}{|k|}))\big)
-
\beta_0\big(2r_1^{-1}(
\inn{\tfrac{-k}{|k|}}{x(\a)}-\rho^*(\tfrac{-k}{|k|}))\big)\Big)
\end{multline*}
and

\begin{multline*}
\Phi_{n}^{\pm}(k,\alpha)=
\frac{\inn{k}{x_\perp'(\a)}}{|k|^2}\quad\times\\
\beta_n\big(R(
\inn{\pm k}{x(\a)}-\rho^*(\pm k)\big)
\beta_0\big(2r_1^{-1}
(\inn{\tfrac{\pm k}{|k|}}{x(\a)}-\rho^*(\tfrac{\pm k}{|k|}))\big).
\end{multline*}

The cutoff function $\Phi_{n}^{+}(k,\cdot)$
 localizes to those points $P$ on the boundary
for which  the
distance of  $P$ to the supporting line
$\{x:\inn{k}{x}=\rho^*(k)\}$ is small and
$\approx 2^n(R|k|)^{-1}$ (or
$\lc (R|k|)^{-1}$ if $n=0$).
Also  $\Phi_{n}^-(k,\cdot)$ gives
a localization in terms of the distance to the supporting line
$\{x:\inn{-k}{x}=\rho^*(-k)\}$.
The factors
$
\beta_0(2r_1^{-1}(
\inn{\tfrac{\pm k}{|k|}}{x(\a)}-\rho^*(\tfrac{\pm k}{|k|})))$
are included in this definition to make sure
 that the supports of $\Phi_{n}^{+}$ and $\Phi_{n}^{-}$ are disjoint.

Note that
$$\Psi(k,\alpha)+\sum_{n=0}^\infty
\Phi_{n}^{+}(k,\alpha)
+\sum_{n=0}^\infty \Phi_{n}^{-}(k,\alpha)=
\frac{\inn{k}{x_\perp'(\a)}}{|k|^2} $$
and also that $\Phi_{n}^{\pm}(k, \alpha)=0$ if $2^n\ge |k|R$.

Define (for fixed $R$ and $h$)
\begin{eqnarray}\label{Inpmk}
I_n^\pm(k,t)&=&2\pi t\,\eta_{R,h}(t)
 \int  \Phi_n^\pm(k,\alpha) e^{- 2\pi \ic \inn{x(\a)}{ tk}} d\a
\\\label{IIk}
II(k,t)&=&2\pi t\,\eta_{R,h}(t)
 \int  \Psi(k,\alpha) e^{-2\pi \ic \inn{x(\a)}{ tk}} d\a
\end{eqnarray}
and
\begin{eqnarray}\label{Inpm}
I^\pm_n(t)&=&\sum_{\frac{2^{n}}{R}<|k|\le R^2}
 \widehat\zeta(2\pi k/R)
I_n^\pm(k,t)
\\
\label{II}
II(t)&=&\sum_{k\neq 0}\widehat \zeta(2\pi  k/R) II(k,t),
\end{eqnarray}
and set
\begin{equation}\label{Gnpm}
G_{n}^\pm(R,h)=\Big(\int |I^\pm_n(t)|^2 dt\Big)^{1/2}.
\end{equation}

Using the decay of $\widehat \zeta$  we see that

\begin{equation}\label{boundsbypieces}
G(R,h)\le \sum_\pm\sum_{n=1}^\infty G_n^\pm(R,h)+
\Big(\int|II(t)|^2 dt\Big)^{1/2} +C R^{-10}.
\end{equation}

\medskip

\noi {\bf Pointwise bounds via van der Corput's lemma.}
We start with a simple pointwise estimate for the pieces
$I^n_\pm$ and $II$ which just
 relies on  van der Corput's lemma for oscillatory integrals (see \cite {S}, p. 334).
We set
\begin{equation} \label{omega}
\om_R(\xi)=(1+|\xi|/R)^{-N}
\end{equation}

\begin{lem}\label{lem:pointwise}
For $n\ge 0$  we have
$$
I_n^\pm(t)\le C 2^{-n}
R|\eta_{R,h}(t)|
\sum_{|k|> 2^n/R}  \omega_R(k) |k|^{-1} \mu(\kunit, \tfrac{2^n}{|k|R})
$$
and
$$II(t)\le C|\eta_{R,h}(t)| \log R.$$
Let
\begin{equation}\label{Gamman}
\Gamma_n^{\pm}(R,h):=
R \sum_{|k|> 2^n/R}  \omega_R(k) |k|^{-1} \mu(\kunit, \tfrac{2^n}{|k|R}).
\end{equation}
Then
\begin{equation}\label{squareptw-a}
G_n^{\pm}(R,h)
\le C\Gamma_n^\pm (R,h);
\end{equation}
moreover
\begin{equation}\label{squareptw-b}
\Big(\int|II(t)|^2 dt\Big)^{1/2}\le C \log R.
\end{equation}
\end{lem}

\medskip
\noi{\bf Proof.}
We write down the argument for $I_n^+$ as  the estimate for $I_n^-$ is
analogous.
The estimate for $n=0$  is immediate if we observe that
length of
the support of $\Phi_0^+(k,\cdot)$ is $\le \mu(k/|k|,(|k|R)^{-1})$.

Fix $\theta\in S^1$ and choose   $\alpha_\theta$ so that
$\inn{\theta}{x(\alpha_\theta)}=\rho^*(\theta)$ and thus also $\fn(\alpha_\theta)=\theta$.
We first observe that if  $\inn{\theta}{x(\a)}-\rho^*(\theta)>\delta$ then
$\tan\measuredangle \big(
\fn(\a),\fn(\a_\theta)\big)\ge \delta/\mu(\theta,\delta) $.

We use this with $\delta=2^n(|k|R)^{-1}$
 to get a lower bound for the derivative of the phase function
in the support of $\Phi_n^*(k,\cdot)$.
This implies that for $t\in \supp \eta_{R,h}$
\begin{equation}\label{lowerbdfirstder}
|\inn{x'(\a)}{ tk}|\ge 2^n \big(\mu(\kunit, \tfrac{2^n}{|k|R})\big)^{-1}
\quad\text{ if }\a\in \supp \Phi_n^*(k,\cdot).
\end{equation}
 and this derivative is monotone in $\a$. Moreover
$$\|\Phi_n^+(k,\cdot)\|_\infty+ \|\partial_\a\Phi_n^+(k,\cdot)\|_1
\lc |k|^{-1};$$
here we use our {\it a priori} assumption on the integrability of the second
derivatives of $\gamma$.

Consequently, by van der Corput's lemma,
we obtain
$$|I_n^+(k,t)|\lc t|\eta_{R,h}(t)|\omega_R(k)
 |k|^{-1}| 2^{-n} \mu(\kunit, \tfrac{2^n}{|k|R})$$
which yields the asserted bound for $I_n^+(t)$.

Similarly,
 $|\partial_\alpha\inn { t k}{x(\alpha)}|\ge c |k|R$  if  $\alpha\in \supp \Psi(k,\cdot)$ and
$\|\Psi(k,\cdot)\|_\infty$ and $ \|\partial_\a\Psi(k,\cdot)\|_1 $ are
 $O(|k|^{-1})$. Thus
$$|II(k,t)|\lc |\eta_{R,h}(t)|  |k|^{-2}(1+|k|/R)^{-N}\omega_R(k)$$
and summing in $k$ yields the asserted bound for $II(t)$. The bounds for the square functions are immediate from the pointwise estimates.\qed

\medskip

\noi{\bf Square function estimates.} We shall need to improve on the pointwise
estimates for $G_n^\pm(R,h)$
 in Lemma \ref{lem:pointwise} which will  only be  useful for  $2^n\ge R$.

We apply Plancherel's theorem with respect to the $t$-variable and obtain
\begin{align*}
G_{n}^\pm&(R,h)^2=2\pi \int |\widehat{I^\pm_n}(\la)|^2 d\la
\\
&=2\pi
\sum_{\frac{ 2^{n}}R<|k|\le R^2}
\sum_{\frac{ 2^{n}}R<|\ell|\le R^2}
\widehat \zeta(2\pi k/R)  \widehat\zeta(2\pi \ell/R)
\int
\widehat{I^\pm_n}(k,\la)
\overline{
\widehat{I^\pm_n}(\ell,\la)}
d\la
\end{align*}
where
\begin{equation}\label{Infourier}
\widehat{I^\pm_n}(k,\la)=
2\pi
\iint t\, \eta_{R,h}(t)
  \Phi_n^\pm(k,\alpha) e^{-\ic t(\la+ \inn{x(\a)}{2\pi k})} d\a dt.
\end{equation}

The crucial estimate is

\begin{lem}\label{lem:secondlemma}
Suppose that $k\in \bbZ^2$ and  $|k|>2^{n}/R$. Then the following inequalities hold.

If $2^n\le R/h$ then
\begin{equation}\label{smalln}
\big|\widehat{I^\pm_n}(k,\la)\big|\lc
\frac{R h^{1/2} |k|^{-1}
 2^{-n}\mu(\pm \kunit,\tfrac{ 2^n}{|k|R})}
{(1+h|\la\pm \rho^*(\pm 2\pi k)|)^{2N}}
\end{equation}
and if $2^n\ge R/h$, then
\begin{equation}\label{largen}
\big|\widehat{I^\pm_n}(k,\la)\big|\lc
\frac{
R h^{1/2} |k|^{-1}
 2^{-n}\mu(\pm \kunit,\tfrac{ 2^n}{|k|R})}
{(1+R2^{-n}|\la\pm\rho^*(\pm2\pi k)|)^{2N}}
\end{equation}

\end{lem}

\medskip

\noi{\bf Proof.} We prove the estimate for $I_n^+(k,\la)$;
the estimate for $I_n^-(k,\la)$ is analogous.

We first consider the case $n=0$. Interchange the order of
integration in (\ref{Infourier}) and perform $2N$
 integrations by parts
 with respect to $t$.  This, together with the estimates
$\Phi_0^+(k,\alpha)=O(|k|^{-1})$, $\eta_{R,h}(t)=O(h^{-1/2})$ yields
\[
|\widehat{I^+_n}(k,\la)|\lc R |h|^{-1/2}|k|^{-1}
\iint\limits_{t\in {\supp}\eta_{R,h} \atop
\alpha\in { \supp}   \Phi_0^+(k,\cdot)}
(1+ h |\la+ \inn{x(\a)}{2\pi k}|)^{-2N} dt d\a.
\]
By definition we have
$|\inn{x(\a)}{k}-\rho^*(k)|\lc R^{-1}$ for $\alpha\in \supp \Phi_0^+(k,\cdot)$. Since $R\ge h$
this implies
\begin{equation}\label{pf22-1}
(1+ h |\la+ \inn{x(\a)}{2\pi k}|)\approx
(1+ h |\la+ \rho^*(2\pi k)|).
\end{equation}
Also observe that the length of the support of $\eta_{R,h}$ is $O(h)$
and that the length of
the support of $\Phi_0^+(k,\cdot)$ is $\le \mu(k/|k|,(|k|R)^{-1})$.
Thus
\begin{eqnarray*}
|\widehat{I^+_0}(k,\la)| & \lc & R |h|^{-1/2}|k|^{-1}
\iint\limits_{ t\in{\supp}\eta_{R,h} \atop
\alpha\in {\supp}   \Phi_0^+(k,\cdot) }
(1+ h |\la+ \rho^*(2\pi k)|)^{-2N} dt d\a\\
&\lc &
R |h|^{1/2}|k|^{-1}\mu(\kunit,\tfrac{1}{|k|R})
(1+ h |\la+ \rho^*(2\pi k)|)^{-2N}
\end{eqnarray*}
which is the asserted estimate for $n=0$.

We now suppose that $n\ge 1$, and begin by performing an integration by parts
with respect to $\alpha$ in (\ref{Infourier}).
Observe that $\inn{x'(\a)}{k}\neq 0$ if $\alpha\in \supp \Phi_n^+(k,\cdot)$.
We obtain
$$I_n^+(k,\la)= F_{n,1}(k,\la)+ F_{n,2}(k,\la)$$ where

\begin{multline}\label{pf22-2}
 F_{n,1}(k,\la)=\\
2\pi
\iint  \eta_{R,h}(t)
  \Phi_n^+(k,\alpha)
\frac{\partial}{\partial \a} \Big (\frac{1}{\ic\inn{x'(\a)}{2\pi k}}\Big)
e^{-\ic t(\la+ \inn{x(\a)}{2\pi k})} d\a dt,
\end{multline}
and
\begin{equation}\label{pf22-3}
 F_{n,2}(k,\la)\\=
2\pi
\iint  \eta_{R,h}(t)
\frac{\partial   \Phi_n^+ }{\partial \a} (k,\alpha)
\frac{e^{-\ic t(\la+ \inn{x(\a)}{2\pi k})} }{\ic\inn{x'(\a)}{2\pi k}}
d\a dt
\end{equation}

As above we interchange the order of integration and integrate by parts in $t$.
This yields the estimate
\begin{multline*}
|F_{n,1}(k,\la)|\lc
R |h|^{-1/2}|k|^{-1} \\ \times
\iint\limits_{t\in{ \supp}\eta_{R,h} \atop
\alpha\in { \supp}   \Phi_n^+(k,\cdot)}
(1+ h |\la+ \inn{x(\a)}{2\pi k}|)^{-2N}
\Big|\frac{\partial}{\partial \a} \Big (\frac{1}{\ic\inn{x'(\a)}{2\pi k}}\Big)\Big|d\alpha dt.
\end{multline*}
If $2^n\le R/h$ then (\ref{pf22-1}) is still valid if $\alpha\in \supp \Phi_{n}^+(k,\cdot)$.
Moreover we claim that
\begin{equation}\label{pf22-4}
\int\limits_{{{ \supp}}   \Phi_n^+(k,\cdot)}
\Big|\frac{\partial}{\partial \a} \Big (\frac{1}{\inn{x'(\a)}{2\pi k}}\Big)\Big|d\alpha \lc
R 2^{-n}\mu(\kunit, \tfrac{2^n}{|k|R}).
\end{equation}

To see this we choose $\alpha_k$ so that $\inn{x(\alpha_k)}{k}=\rho^*(k)$
(this choice may not be unique).  The support of $\Phi_n^+(k,\cdot)$ consists of two
connected intervals (on ${\mathbb R}/{L}{\mathbb Z}$ with  $L=|\partial \Omega|$)
 and on each of these
the function
 $\alpha\to \inn{x'(\a)}{k}$ is monotone; moreover this function vanishes for
$\alpha=\alpha_k$. Thus
\begin{eqnarray}
2^{n-1}R^{-1} &\le & \big|\inn{x(\a)}{k}-\rho^*(k)\big| \nonumber \\
& = & \Big|(\alpha-\alpha_k)\int_0^1\inn{x'(\alpha_k+\tau(\a-\a_k))}{k}
d\tau\Big|
\label{pf22-5}\\
&\le & |\alpha-\alpha_k| |\inn {x'(\a)}{k}|. \nonumber
\end{eqnarray}

Note that $|\alpha-\alpha_k|\lc \mu(k/|k|, 2^n/(|k|R))$ for $\alpha\in\supp
\Phi_n^+(k,\cdot)$.
Since
$\partial_\alpha
((\inn{x'(\a)}{2\pi k})^{-1})$ is single-signed on the two components of
$\supp \Phi_n^+(k,\cdot)$  we can apply the fundamental
theorem of calculus on these intervals
and we see that the left hand side of (\ref{pf22-4}) is bounded by
$4\sup|\inn{x'(a)}{2\pi k}|^{-1})$ where the supremum  is
taken over all $\alpha \in
\supp(\Phi_n^+(k,\cdot))$. But by (\ref{pf22-5}) this bound is
$O(R2^{-n}\mu(k/|k|, 2^n/(|k|R)))$.

 Combining (\ref{pf22-1}) and (\ref{pf22-4}) we obtain for
$2^n\le R/h$
\begin{multline}
\label{pf22-6}
|F_{n,1}(k,\la)| \\
\lc  |h|^{1/2}|k|^{-1}
(1+ h |\la+ \rho^*(2\pi k)|)^{-2N}
\int\limits_{
 \supp   \Phi_n^+(k,\cdot)}\Big|
\frac{\partial}{\partial \a} \Big (\frac{1}{\inn{x'(\a)}{2\pi k}}\Big)\Big|d\alpha \\
\lc 2^{-n} R |h|^{1/2}|k|^{-1}
(1+ h |\la+ \rho^*(2\pi k)|)^{-2N}
\mu(\kunit, \tfrac{2^n}{|k|R})
\end{multline}
which is the estimate we were aiming for.
Next we consider the term $F_{n,2}(k,\la)$ and arguing as above we see that
\begin{equation} \label{pf22-7}
|F_{n,2}(k,\la)|\lc  h^{1/2}
(1+ h |\la+ \rho^*(2\pi k)|)^{-2N}
\int\Big|\frac{\partial   \Phi_n^+ }{\partial \a} (k,\alpha)\Big|
\Big|\frac{1}{\inn{x'(\a)}{2\pi k}}\Big|
d\a
\end{equation}
and
\begin{eqnarray*}
& &\int\Big|\frac{\partial   \Phi_n^+ }{\partial \a} (k,\alpha)\Big|
\Big|\frac{1}{\inn{x'(\a)}{2\pi k}}\Big|
d\a
\\
& & \quad \lc \int_{\supp \Phi_n^+(k,\cdot)}
 \frac{|\inn{k}{x_\perp''(\a)}|}{|k|^2}
\Big|\frac{1}{\inn{x'(\a)}{2\pi k}}\Big|d\a
\\ &&\quad +
\frac 1{|k|}
\int_{\supp \Phi_n^+(k,\cdot)}
\big|\partial_\alpha\big(\beta(R2^{-n}(\inn{x'(\a)}{k}-\rho^*(2\pi k)))\big)\big| d\a
\\ && \quad+
\frac 1{|k|}
\int_{\supp \Phi_n^+(k,\cdot)}
\big|\partial_\alpha
\big(\beta_0(2r_1^{-1}(
\inn{\tfrac{k}{|k|}}{x(\a)}-\rho^*(\tfrac{k}{|k|})))\big)
\big| d\a
\\&& := A_1(k)+A_2(k)+A_3(k).
\end{eqnarray*}
We now use  that by (\ref{pf22-4})
$$|\inn{x'(a)}{2\pi k}|^{-1}\lc R 2^{-n}\mu(k/|k|,2^n/(|k|R))$$
 on the support of
$\Phi_n^+(k,\cdot)$.
Since $\inn{k}{x_\perp''(\a)}$ is single-signed on the components we see that
$$
A_1(k)\lc R 2^{-n} |k|^{-1} \mu(\kunit, \tfrac{2^n}{|k|R}).
$$
Next
$$A_2(k)\lc |k|^{-1} R2^{-n} \meas(\supp  \Phi_n^+(k,\cdot))
\lc |k|^{-1} R 2^{-n}\mu(\kunit, \tfrac{2^n}{|k|R}).
$$
Finally, on the support of the  derivative of the term  $\beta_0(\dots)$ we have the better bound
$|\inn{x'(\a)}{k}|^{-1}=O( (|k|)^{-1})$ so that
\begin{align*}
A_3(k)&\lc |k|^{-2} \meas(\supp  \Phi_n^+(k,\cdot))
\\&\lc  |k|^{-2} \mu(\kunit, \tfrac{2^n}{|k|R})
\lc |k|^{-1} R 2^{-n}\mu(\kunit, \tfrac{2^n}{|k|R})
\end{align*}
in view of our restriction $k\ge 2^n/R$. If we use
these estimates in (\ref{pf22-7})
  then we obtain the desired estimate for
$F_{n,2}(k,\la)$, at least for the case
$2^n\le R/h$.

The estimates for
$F_{n,1}(k,\la)$ and $F_{n,2}(k,\la)$ in the case
$2^n> R/h$ are derived analogously. The only difference is that
(\ref{pf22-1}) does not hold in all
of the support of $\Phi^\pm_n(k,\cdot)$. However we still
 have $|\rho^*(2\pi k)-\inn{x'(\a)}{2\pi k}|\lc 2^n/R$
in this set so that (\ref{pf22-1}) is now replaced by
\begin{multline} \label{pf22-8}
(1+ h |\la+ \inn{x(\a)}{2\pi k}|)\\
\lc
(1+ R2^{-n} |\la+ \inn{x(\a)}{2\pi k}|)\approx
(1+ R2^{-n} |\la+ \rho^*(2\pi k)|)
\end{multline}
and the remainder of the above arguments
applies without change to yield the inequalities in (\ref{boundsbypieces})
\qed


\begin{lem}\label{thirdlemma}
Let
\begin{multline}\label{fBn}
\fB^n_\pm(R,h):=
 2^{-n} R \, \times\, \\
 \Big(\sum_{k\in {\mathbb Z}^2 \atop \frac{2^n}{R}< |k|\le R^2}
\sum_{\ell\in {\mathbb Z}^2 \atop
 \frac{2^n}{R}< |\ell|\le R^2}
\frac{\om_R(k)\om_R(\ell)}{|k||\ell|}
\frac{\mu(\kunit, \tfrac {2^n}{|k|R})\,
\mu(\lunit, \tfrac {2^n}{|\ell|R})}
{(1+h|\rho^*(\pm k)-\rho^*(\pm\ell)|)^N}
\Big)^{1/2},
\end{multline}
and

\begin{multline}\label{fBntilde}
\widetilde{\fB}^n_\pm(R,h):=
 2^{-n} R \, \times\, \\
\quad \Big(\sum_{k\in {\mathbb Z}^2 \atop \frac{2^n}{R}< |k|\le R^2}
\sum_{\ell\in {\mathbb Z}^2 \atop
 \frac{2^n}{R}< |\ell|\le R^2}
\frac{\om_R(k)\om_R(\ell)}{|k||\ell|}
\frac{\mu(\kunit, \tfrac {2^n}{|k|R})\,
\mu(\lunit, \tfrac {2^n}{|\ell|R})}
{(1+R2^{-n}|\rho^*(\pm k)-\rho^*(\pm\ell)|)^N}
\Big)^{1/2}.
\end{multline}

Then  for $1\le h\le R$, $R\ge 2$,
\begin{align}
\label{Gnfirst}
G_n^\pm (R,h) &\le C \fB_n^\pm(R,h) \quad \text{ if } 2^n\le R/h,
\\
\label{Gnsecond}
G_n^\pm (R,h) &\le C \widetilde{\fB}_n^\pm(R,h) \quad \text{ if } 2^n> R/h.
\end{align}
\end{lem}

\medskip

\noi{\bf Proof.}
 We observe that $|\zeta(k/R)|\le \om_R(k)$ and use the elementary
convolution inequality
\begin{equation}\label{convol}
\int
{(1+a|A+\la)|)^{-2N}}
{(1+a|B+\la)|)^{-2N}}d\la \lc a^{-1} (1+a|A-B|)^{-N}
\end{equation}
We use (\ref{convol}) for $A=\rho^*(\pm 2\pi k)$,
$B=\rho^*(\pm 2\pi \ell)$ and
$a=h$ if $2^n \le R/h$ and  $a=R2^{-n}$ if
$ 2^n > R/h$.
The estimates (\ref{Gnfirst}), (\ref{Gnsecond}) are   now  straightforward from
(\ref{smalln}), (\ref{largen})  and (\ref{convol}).\qed

We shall now combine the previous lemmata to  state the definitive
 result of this section; here we abandon the {\it a priori} assumption that the outer unit normals
are absolutely continous functions.

\begin{prop}\label{prop:2.4} Let $\Omega$ be a convex domain
containing the origin in its
interior, and suppose that $r_1< 1<r_2$ where  $r_1,r_2$
 are the radii of inscribed and circumscribed circles centered at the origin.
Let
$\fB_n^\pm$,
$\widetilde{\fB}_n^\pm$ and $\Gamma_n^\pm$ be as in (\ref{fBn}), (\ref{fBntilde}) and (\ref{Gamman}).

 There exists a constant
 $C$, depending only on $r_1$, $r_2$ and $N$, so that
 for $1\le h\le R$, $R\ge 2$ we have the estimate
\begin{multline}\label{main-sect2}
\cG_\Om(R,h)\le C\sum_\pm\Big[
  \sum_{2^n\le R/h} \fB_n^\pm(R,h)+
  \sum_{R/h<2^n\le R} \widetilde{\fB}_n^\pm(R,h)\\+
  \sum_{2^n>R} \Gamma_n^\pm(R,h)\Big]+C [\log R+(R/h)^{1/2}].
\end{multline}
\end{prop}

\medskip

\noi{\bf Proof.} Under our previous  {\it a priori} assumption on the
 boundary of $\Omega$ this statement
 follows by simply putting together the estimates
(\ref{Emollified}), (\ref{boundsbypieces}),
(\ref{squareptw-a}),
(\ref{squareptw-b}),
(\ref{Gnfirst}), and
(\ref{Gnsecond}). We note that all bounds just depend on $r_1\le 1\le r_2$
and $N$,
and that the $L^1$ bound  for the second derivatives does not enter
in the result.

In the general case  we note that
there is a sequence of convex domains $\Omega_j$ which contain
the origin, such that
 $\Omega_j\subset \Omega_{j+1}\subset \Omega$ and
$\cup_j\Omega_j=\Omega$ and $\Omega_j$ has smooth  boundary.
Moreover, let
$\mu_j(\theta,\delta)$ for
fixed   $\theta\in S^1$ and $\delta>0$  be the quantity (\ref{eq:1.6}) but
associated to
the domain $\Om_j$. Then $\mu_j(\theta,\delta)$  converges to the corresponding
quantity associated to $\Omega$, $\mu(\theta,\delta)$.
 Moreover the square functions  defined by the
smoothed errors $E_{1/R}$  associated to
$\Omega_j$ converge to the corresponding expression associated to $\Omega$ and
the same statement applies to the expressions
$\fB_n^\pm(R,h)$, $\widetilde{\fB_n}^\pm(R,h)$, $\Gamma_n^\pm(R,h)$.
 For one explicit construction of the approximation see the proof of
 Lemma 2.2 in \cite{SZ}.
The  constant $C$  in the statement of Lemmata \ref{lem:pointwise} and
\ref{thirdlemma}  can be
chosen uniformly in $j$ and the assertion follows.\qed

\bigskip

\section{Estimates for the case of nonzero curvature}
\label{sectionnonzerocurvature}

In this section we estimate the various quantities 
in \eqref{main-sect2} of Proposition \ref{prop:2.4},
in the case of nonvanishing curvature
and rough boundary; that is, we assume inequality (\ref{eq:1.7}).
Proposition \ref{prop:2.4} reduces matters to estimates for
lattice points in  thin annuli
\begin{equation}\label{annuli} \cA^\pm (r,h):=
\{\xi\in \bbR^2: |\rho^*(\pm\xi)-r|\le h^{-1}\};
\end{equation}
here $r\ge 2$ and $h\ge 1$. Let
\begin{equation}\label{Srh}
S^\pm (r,h):= \card (\cA^\pm(r,h)\cap \bbZ^2),
\end{equation}
the number of lattice points in the annulus
$ \cA^\pm (r,h)$.

We also need to consider  for $\delta\le r^{-1}$ the weighted sum
\begin{equation}
\label{weightedsum} \fS^\pm(r,\delta,h)=
\sum_{\ell\in \bbZ^2:\atop
|\rho^*(\pm\ell)-r|\le h^{-1}}
\frac{ \mu(\ell/|\ell|, \delta)}
{\sqrt \delta}.
\end{equation}

Notice that assuming
(\ref{eq:1.7}) we have
\begin{equation}
\label{trivialestimate}
\fS^\pm(r,\delta,h) \lc S^\pm (r,h);
\end{equation}
however
$\fS^\pm(r,\delta,h)$ {\it  could be much smaller} than $S^\pm(r,h)$.
Indeed consider the case that
 $S^\pm(r,h)$ contains a  long line segment;  then the boundary
$\Omega^*$ becomes  nearly flat and by duality the curvature of
$\Omega$  at the corresponding points gets large causing
$\mu(\tfrac{\ell}{|\ell|},\delta)$ to be smaller than $\delta^{1/2}$ for many
$\ell$ on the line segment.
This phenomenon  will be  exploited in the proof of Theorem 1.2.


\begin{prop}\label{prop:3.1}
Suppose that assumption \eqref{eq:1.7} is satisfied and that
 $R\ge 2$ and $1\le h\le R/\log^2 R$.
Then
\begin{multline}
\label{mainbound}
\cG_\Omega(R,h)
\le C R^{1/2} \\+CR^{1/2}  \sum_\pm\sum_{0\le 2^n \le \log R}2^{-n/2}\Big(
\sum_{1\le l\le \log R}
2^{-l}\sup_{2^{l-1}\le r\le 2^{l}} \fS^\pm\big(r, \tfrac{2^n}{rR}, h\big)
 \Big)^{1/2}.
 \end{multline}
\end{prop}

\noi{\bf Proof.} The proof relies on rather
straightforward calculations 
used to bound the expressions on the right hand side of  \eqref{main-sect2}.
We  shall employ the bounds 
$$\mu(\kunit,\tfrac{2^n}{|k|R})=O(2^{n/2} (R|k|)^{-1/2})$$
(by \eqref{eq:1.7}) and  the definition of $\fS^\pm$.

We now proceed to estimate the various terms in \eqref{main-sect2}.
First, in order  to bound the terms $\Gamma_n^\pm(R,h)$, in the range  $2^n\ge R$ 
we get from \eqref{Gamman} 
$$
\Gamma_n^\pm(R,h)
\lc 2^{-n/2} R^{1/2}\sum_{|k|>0} |k|^{-3/2}\omega_R(k)
\lc R2^{-n};
$$
Thus the contribution of
$\sum_{2^n\ge R}
\Gamma_n^\pm(R,h)$ is certainly covered  by the first term on the right hand side of 
\eqref{mainbound}.

Next we consider the terms
$\widetilde\fB_n^\pm(R,h)$ which were defined in \eqref{fBntilde}.
Now  $R/h<2^n\le R$, and
 in view of our assumption
$h\le R(\log R)^{-2}$ this certainly implies $2^n\ge \log R$.
%
We compute
from (\ref{fBntilde}),
using  the assumption $R2^{-n}\ge 1$,
\begin{align*}\widetilde {\fB}_n^\pm(R,h)&\lc
 2^{-n/2} R^{1/2}
\quad \Big(\sum_{k\in {\mathbb Z}^2 \atop \frac{2^n}{R}< |k|\le R^2}
\sum_{\ell\in {\mathbb Z}^2 \atop
 \frac{2^n}{R}< |\ell|\le R^2}
\frac{\om_R(k)\om_R(\ell)(|k||\ell|)^{-3/2}}
{(1+|\rho^*(\pm k)-\rho^*(\pm\ell)|)^N}
\Big)^{1/2}
\\&\lc h^{1/2}\log^{1/2} R, \qquad \text{if }\, R/h\le 2^n\le R.
\end{align*}
Therefore,
$$
\sum_{R/h\le 2^n\le R}\widetilde {\fB}_n^\pm(R,h)\lc
(h\log(2+h^{-1})\log R)^{1/2} \lc R^{1/2}
$$
since we assume $h\le R/\log^2R$. Thus these contributions are again subsumed 
under the first bound in \eqref{mainbound}.

Next 
we bound
${\fB}_n^\pm(R,h)$, defined in \eqref{fBn}, now considered  for the range 
$2^n\le R/h$. 
The argument above for $\widetilde {\fB}_n^\pm(R,h)$  also applies to
${\fB}_n^\pm(R,h)$
and one gets ${\fB}_n^\pm(R,h)=O(2^{-n/2} R^{1/2} (\log R)^{1/2})$
if
$R/h\le 2^n\le R$.
Thus $$
\sum_{\log R \le 2^n\le R/h}{\fB}_n^\pm(R,h)\lc R^{1/2}.
$$
Note that if instead  we summed  over the range
 $1\le 2^n\le R/h$ then  we would only
get the weaker bound $O(R^{1/2}(\log R)^{1/2})$.

Finally we have to bound
$\sum_{2^n\le \log R}{\fB}_n^\pm(R,h)$. 
To 
this end we observe  that
 the sum of the  contributions of the terms in (\ref{fBn})
which involve either $|k|\ge R$ or $|\ell|\ge R$ or
$|\rho*(\pm k)-\rho^*(\pm \ell)|\ge \sqrt{\rho^*(\pm k)}$
is easily recognized to be $O(R^{1/2})$.

Now let $E^\pm(k)$ denote the  set of all
$\ell \in \mathbb Z^2$ which also satisfy
$0<|\ell|\le R$ and
$|\rho^*(\pm k)-\rho^*(\pm \ell)|<\sqrt{\rho^*(\pm k)}$. We use  the bound
$\mu( \kunit, \tfrac{2^n}{|k|R})\lc 2^{n/2}|k|^{-1/2}R^{-1/2}$ and estimate
\begin{multline*}{\fB}_n^\pm(R,h)\le C_1  R^{1/2}\quad +
\\   C_2
 2^{-n3/4} R^{3/4}
\quad \Big(\sum_{k\in {\mathbb Z}^2 \atop 0< |k|\le R}|k|^{-5/2}
\sum_{ \ell\in E^\pm (k,n)}
\frac{\mu(\lunit, \tfrac{2^n}{|\ell|R})}
{(1+h|\rho^*(\pm k)-\rho^*(\pm\ell)|)^{N}}
\Big)^{1/2}.
\end{multline*}
By the property $\mu(\theta, A\delta)\le
 C_A\mu(\theta,\delta)$ we obtain

\begin{align*}
&\sum_{ \ell\in E^\pm (k)}
\frac{\mu(\lunit, \tfrac{2^n}{|\ell|R})}
{(1+h|\rho^*(\pm k)-\rho^*(\pm\ell)|)^{N}}
\\&\lc
\sum_{|m|\le C|k|^{1/2}}  (1+|m|)^{-N}
\Big(\frac {2^n}{R|k|}\Big)^{1/2}
\fS^\pm \big(\rho^*(\pm k)+m,
\tfrac{2^{n}}{|k|R}, h\big)
\end{align*}
and thus
\begin{align*} &2^{-3n/4} R^{3/4}
\quad \Big(\sum_{k\in {\mathbb Z}^2 \atop 0< |k|\le R}|k|^{-5/2}
\sum_{ \ell\in E^\pm (k)}
\frac{\mu(\lunit, \tfrac{2^n}{|\ell|R})}
{(1+h|\rho^*(\pm k)-\rho^*(\pm\ell)|)^{N}}
\Big)^{1/2}
\\
&\lc 2^{-n/2} R^{1/2}
\quad \Big(\sum_{k\in {\mathbb Z}^2 \atop 0< |k|\le R}|k|^{-3}
\sup_{\frac{\rho^*(\pm k)}2\le r\le 2\rho^*(\pm k)}
\fS^\pm(r,\tfrac{2^{n+1}}{rR},h) \Big)^{1/2}
\end{align*}
which is in turn bounded by a constant times
$$
2^{-n/2} R^{1/2}
\quad \Big(\sum_{1\le l\le 2+
\log R} 2^{-l} \sup_{2^{l-1}\le r\le 2^l}
\fS^\pm(r,\tfrac{2^{n}}{rR},h)\Big)^{1/2}.
$$\qed


We now state the crucial propositions needed in the proof of Theorems
\ref {thm:1.1} and
\ref {thm:1.2}.
The mild regularity assumption in Theorem  \ref {thm:1.1} gives us a
favorable  estimate for  $S^\pm(r,h)$ which will be proved in
\S\ref{sectionlogregularity}.

\begin{prop}\label{prop:3.2} Let $\Omega$
be as in the statement of Theorem 1.1, i.e. with
$\kappa\in L\log^{2+\epsilon} L$ and $\kappa$ bounded below.
Assume that  $2\le r$ and $1\le h$. Then
\begin{equation}
S^\pm(r,h)\lc r\big[ h^{-1}+ \log^{-1-\eps/2}(2+r)\big].
\end{equation}
\end{prop}

If we only make the assumption that the curvature is bounded below
(in the sense of
\eqref{eq:1.7}) then there is no nontrivial  pointwise bound
for
$S^\pm(r,h)$, but we still have a favorable bound for
the weighted sums
$\fS(r, (Rr)^{-1}, h)$ provided that $h\gc R^{1/2}$,
\cf. the comment following \eqref{trivialestimate}.
The bound for these weighted  sums
is   more difficult than Proposition \ref{prop:3.2} and
the combinatorial proof of the following proposition  will be given in
\S\ref{sectiondiscrete}-\ref{proofprop33}.
It may fail for $h<R^{1/2}$,
see \S\ref{sharpness}.

\begin{prop}\label{prop:3.3}
Let $\Omega$
be as in the statement of Theorem 1.2, i.e. with
$\kappa$ bounded below.
Assume that $R\ge 10$, $10\le r\le R$ and $R^{1/2}\le h\le  R$.

Then
\begin{equation}\label{equationinprop3.3}
\sqrt{Rr} \,
\fS(r, \tfrac 1{Rr}, h)
\lc r^{17/18}\,.
\end{equation}
\end{prop}

%

%

We finish this section by showing how the above  propositions imply the
results stated in the introduction.

\medskip

\noi{\bf Proof of Theorem 1.3}. For this result we just use the trivial estimate $S^\pm (r,h)=O(r)$ if $r\ge 1$, $|h|\ge 1$. The bound
$\cG_\Om(R,h)=O((R\log R)^{1/2})$ follows easily from a combination of
Proposition \ref{prop:3.1} and  \eqref{trivialestimate}.
\qed

\medskip

\noi{\bf Proposition 3.2 implies  Theorem 1.1.}
Now we  still use \eqref{trivialestimate} and observe that by
Proposition \ref{prop:3.2}
\begin{multline*}
\sum_{1\le l\le \log R}
2^{-l}\sup_{2^{l-1}\le r\le 2^{l}} \fS\big(r, \tfrac{2^n}{rR}, h\big)
\\
 \lc \sum_{1\le l\le \log R}(h^{-1}+(1+l)^{-1-\eps/2})
 \lc (1+h^{-1}\log R)
\end{multline*}
and thus we obtain the bound $\cG_\Om(R,h)=O(R^{1/2})$ from Proposition
\ref{prop:3.1} if $h\ge \log R$.\qed

\medskip

\noi{\bf Proposition 3.3 implies  Theorem 1.2.}
We argue similarly but now use
the  inequality
\eqref{equationinprop3.3}
 in the application of
 Proposition \ref{prop:3.1}.
Here
\eqref{equationinprop3.3}
is applied with $R$ replaced by
$R2^{-n}$, and since $2^n\le \log R$ this application is certainly valid
for $R^{1/2}\le h\le R/(\log R)^2$.  We obtain
\begin{multline}\label{maintermwithequationinprop3.3}
\sum_{1\le l\le \log (R2^{-n})}
2^{-l}\sup_{2^{l-1}\le r\le 2^{l}} \fS\big(r, \tfrac{2^n}{rR}, h\big)
\\
 \lc \sum_{1\le l\le \log (R2^{-n})} \big( 2^{-l/18}+
(\log(R2^{-n}))^{-1} \big)
\le C;
\end{multline}
moreover for the terms with $\log (R2^{-n})< l\le \log R$ we simply use the trivial bound
$\fS\big(r, \tfrac{2^n}{rR}, h\big)=O(r)$ and get
\begin{equation}\label{maintermwithoutequationinprop3.3}
\sum_{\log (R2^{-n})<l\le \log R}
2^{-l}\sup_{2^{l-1}\le r\le 2^{l}} \fS\big(r, \tfrac{2^n}{rR}, h\big)
\lc n+1.
\end{equation}
We use \eqref{maintermwithequationinprop3.3}   and
\eqref{maintermwithoutequationinprop3.3}  in the application of
Proposition
\ref{prop:3.1}, and in view of the exponential decay in $n$ in
\eqref{mainbound} we
 obtain the bound $\cG_\Om(R,h)=O(R^{1/2})$.\qed

\section{Proof of Theorem \ref{thm:1.4}}
We follow the same setup as in the proof of Theorem
\ref{thm:1.3} and use a
crucial  fact from  \cite{BCI} according to which the maximal function
defined by
\begin{equation*}
\mu^*(\theta)=\sup\{\delta^{-1/2}\mu(\theta, \delta): \delta>0\}
\end{equation*}
belongs to $L^{2,\infty}(S^1);$
{\it i.e.}
\begin{equation*}
\text{meas}(\{\theta\in S^1: \mu^*(\theta)^2>s\})\le C^2/s
\end{equation*}
uniformly in $s$.

We now consider the sets $A_\vartheta\Om$
and denote the quantities in \eqref{Gnpm} associated to $A_\vartheta\Omega$ by
$G_{n}^\pm(R,h, \vth)$ etc. By averaging
it suffices to assume $h=1$. We estimate
$G_{n}^+(R,1, \vth)$.

\par From  Lemma \ref{thirdlemma}  we obtain in the range  $2^n\le R$
\begin{multline*}
G_{n}^+(R,1,\vth)  \lc  2^{-n/2} R^{1/2}  \\
  \times
\Big(\sum_{0<|k|\le R^2 \atop 0< |\ell|\le R^2 }
\mu^*(A_\vth\kunit) \mu^*(A_\vth\lunit)
\frac{|k|^{-3/2}|\ell|^{-3/2}}{
(1+|\rho^*(A_\vth k)-\rho^*(A_\vth\ell)|)^{N}}
\Big)^{1/2}.
\end{multline*}
By symmetry we may restrict the summation to those pairs $(k,\ell)$ for which $\mu^*(\kunit)
\le \mu^*(\lunit)$ and we thus  have the estimate
\begin{multline*}
G_{n}^+(R,1,\vth)\lc 2^{-n/2} R^{1/2}\\
\times
\Big(\sum_{0<|k|\le R^2 \atop 0< |\ell|\le R^2 }
\mu^*(A_\vth\kunit)^2
\frac{|k|^{-3/2}|\ell|^{-3/2}}
{(1+|\rho^*(A_\vth k)-\rho^*(A_\vth\ell)|)^{N}}
\Big)^{1/2}.
\end{multline*}
Now as above it is easy to see that for fixed $k$
$$
\sum_{\ell\neq 0}
|\ell|^{-3/2}
(1+|\rho^*(A_\vth k)-\rho^*(A_\vth \ell)|)^{-N}
\le C |k|^{-1/2}
$$
where $C$ is independent from $\vth$. Thus
\begin{equation}\label{Gn-estimate}
G_{n}^+(R,1,\vth)
\lc 2^{-n/2} R^{1/2}
\sum_{0<|k|\le R^2}
\mu^*(A_\vth\kunit)^2
|k|^{-2}
(1+\tfrac{|k|}R)^{-N};
\end{equation}
moreover
\begin{eqnarray}
& & \sup_{j\ge n}
2^{-j} j^{-(2+\epsilon)} \sup_{2^j\le R\le 2^{j+1}}
G_{n}^+(R,1,\vth)^2
\nonumber  \\
& & \lc
2^{-n}
\sum_{j=n}^\infty
 j^{-(2+\epsilon)}
\sum_{0< |k|\le 2^{2j}}
\mu^*(A_\vth\kunit)^2
|k|^{-2}
(1+2^{-j}|k|)^{-N}.
  \label{Gn-max}
\end{eqnarray}
In order to complete the proof we
have to show that the  expression (\ref{Gn-max})
defines a function in
$L^{1,\infty}([-\pi,\pi])$.
We  apply a well known lemma by Stein and N. Weiss \cite{SW} on adding functions in $L^{1,\infty}$
and the quasi-
norm is bounded by a constant times the square-root of
\begin{eqnarray*}
&  & 2^{-n} \sum_{j=n}^\infty  j^{-(2+\epsilon)}
\sum_{0< k\le 2^{2j} }
|k|^{-2}(\log (1+|k|+j)) \nonumber\\
& & \lc 2^{-n} \sum_{j=n}^\infty  j^{-(1+\epsilon)}\le C_\eps 2^{-n}.
\end{eqnarray*}
Thus
the function
\begin{equation}
\vth\mapsto \sup_{R\ge 2^n} R^{-1/2}(\log(2+ R))^{-1-\eps}
G_{I,n}^+(R,1,\vth )
\label{maxfct}
\end{equation}
belongs to $L^{2,\infty}$ with norm $O(2^{-n/2})$.

For $R<2^n\le R^3$ we argue as in the proof of
Theorem \ref{thm:1.2} and see that
the  estimate (\ref{Gn-estimate}) is replaced by
\begin{equation*}
G_{n}^+(R,1,\vth)
\lc
\sum_{ 0<|k|\le R^2}
\mu^*(A_\vth\kunit)^2
|k|^{-2}
\end{equation*}
and thus
\begin{eqnarray*}
& &  \sup_{j< n}
2^{-j\eps}  \sup_{2^j\le R\le 2^{j+1}}
G_{n}^+(R,1,\vth)^2\\
&  & \quad \lc
\sum_{n/3<j\le n}
2^{-j\eps}
\sum_{2^{n-j}< |k|\le 2^{2j}}
\mu^*(A_\vth\kunit)^2
|k|^{-2}
\end{eqnarray*}
and again this expression as a function of $\vth$ belongs to $L^{1,\infty}$ with quasi-norm
$2^{-\eps n/3}$.
Thus the function
$$\vth\mapsto \sup_{R< 2^n} R^{-\eps/2}
G_{n}^+(R,1,\vth)
$$
belongs to $L^{2,\infty}$ with norm $O(2^{-\eps n/3})$
(which is a better result than for the function
(\ref{maxfct}), as was to be expected).
We may sum in $n$ and get the required assertion for $\vth\mapsto \sum_{n=0}^\infty G_n^+(R,1,\vth)$
and the corresponding assertion involving
$G_n^-(R,1,\vth)$ follows in the same way.
\qed

\section{Bounds for the lattice rest  associated to the polar set --\\
the proof of Proposition \ref{prop:3.2}}
 \label{sectionlogregularity}

We improve
the trivial estimate
$E_{\Om^*}(t)=O(t)$ under the given mild regularity
assumption on $\partial\Omega$.
Because of
$$S^+(r,h)\lc rh^{-1}+ \sup_{t\le r+h^{-1}} E_{\Om^*}(t)$$
(and a similar estimate for $S^{-}$)
 Proposition \ref{prop:3.2} is an immediate consequence of  the
following result.

\begin{prop}\label{prop:5.1}
 Let $\Omega$
be a convex domain with $C^1$ boundary in $\bbR^2$ containing the
origin in its interior, and assume that the
components of the tangent vector are absolutely continuous.
 Suppose also that the  curvature $\kappa$ is uniformly bounded below,
i.e. $|\kappa(x)|\ge a>0$ for almost every $x\in\partial \Om$
and that $\kappa\in L\log^{\ga}\!L(\partial\Omega)$, for some $\ga>0$.
Let $E_{\Om^*}(t)=\cN_{\Om^*}(t)-t^d\volume(\Om^*)$.
Then for $t\ge 2$
\begin{equation}
  |E_{\Om^*}(t)|\le C  t(\log t)^{-\ga/2}\label{eq:5.1}
\end{equation}
\end{prop}

We need the following variant of van der Corput's Lemma.

\begin{lem}\label{lem:5.2} Let $f$ be a $C^1$ function on the interval $[a,b]$ and
assume that $f'$ is absolutely continuous and monotone. Let $\ga>0$ and
suppose that
the function
$t\mapsto (\log (2+ \tfrac 1{|f''(t)|}))^{\ga}$ belongs to $L^{1,\infty}$, with operator (quasi-) norm bounded by $A$.
Then
$$\Big|\int_a^b e^{\ic\la f(t)}\chi(t) dt\Big|\le C(\ga,A)
\big(\|\chi\|_{\infty}+\|\chi'\|_1\big) (\log(2+\la))^{-\ga}.$$
\end{lem}

\noindent {\bf Proof.}
We may assume that $\la \ge 10$.
In view of the monotonicity of $f'$ the set
$I=\{t\in [a,b]:|f'(t)|\le \la^{-1}(\log\la)^\ga\}$ is an interval, $I=[c,d]$.
The set $[a,b]\setminus I$ is a union of at most two intervals
and on each of these  we have $|\la  f'(t)|\ge (\log \la)^\ga$. By
  the standard van der Corput Lemma with first derivatives (\cite{SW})
it follows that
\begin{equation}
\Big|\int_{[a,b]\setminus I} e^{\ic\la f(t)}\chi(t) dt\Big|\le C
\big(\|\chi\|_{\infty}+\|\chi'\|_1\big) (\log\la)^{-\ga}.
\label{eq:5.2}
\end{equation}

To complete the proof we have to show that
\begin{equation}
|I|\lc (\log\la)^{-\ga}.
\label{eq:5.3}
\end{equation}
 Let $E_1=\{ t\in I: |f''(t)|\le (\log\la)^{2\ga}\la^{-1}\}$ and $E_2=I\setminus E_1$.
On $E_1$ we have
$$\log^\ga(2+\tfrac 1{|f''(t)|})\ge \log^\ga (2+\tfrac{\la}{\log^{2\ga} \la})\ge c
\log^\ga(2+\la);
$$
here $c$ depends only on $\gamma$.
Thus by our $L^{1,\infty}$  assumption $|E_1|\lc(\log (2+\la))^{-\ga}$.
By definition of $I$, and the a.e. nonnegativity of $f''$  we also have
$$
2\frac{(\log\la)^{\ga}}{\la}
\ge |f'(d)-f'(c)|
=\Big|\int_I f''(s) ds\Big|\ge
\Big|\int_{E_2} f''(s) ds\Big|\ge |E_2|\frac{(\log\la)^{2\ga}}{\la}
$$
thus $|E_2|\le 2(\log \la)^{-\ga}$.
Thus we have shown (\ref{eq:5.3}) and the proof is complete.
\qed

As a consequence we obtain

\begin{lem}\label{lem:5.3} Let $\Omega$ be as in Proposition
\ref{prop:5.1}.
Then
\begin{equation}
\big|\widehat {\chi_{\Omega^*}}(\xi)\big|\le  C (1+|\xi|)^{-1} \log (2+|\xi|)^{-\ga}
\label{eq:5.4}
\end{equation}
\end{lem}

\noindent{\bf Proof.}
Let $\alpha\mapsto x(\alpha)$ be a parameterization of $\partial \Omega$ with $|x'(\a)|=1$.
A parametrization of $\partial \Omega^*$ is then given by
$\alpha\mapsto \widetilde x(\a)=\inn{x(\a)}{\fn(\a)}^{-1} \fn(\a)$ but this parametrization
is not sufficiently regular.
We compute
\begin{equation}
\widetilde x'= \frac{n_1'n_2-n_1 n_2'}{\inn{x}{\fn}^2}(x_2, -x_1)
\label{eq:5.5}
\end{equation}
and we observe that $n_1'n_2-n_1 n_2'=\ka$. Moreover, if
$r_1<r_2$ are the radii of  inscribed and circumscribed circles centered at the origin then
for $x\in \partial \Om$
\begin{equation}
\inn{x}{\fn}\ge \frac{r_1}{2r_2}|x| >\frac{r_1^2}{2r_2}.
\label{eq:5.6}
\end{equation}
We introduce a new parameter $\tau=\tau(\alpha)=\int_{\alpha_0}^\alpha\frac{ \kappa(\beta)}
{\inn{x(\beta)}{\fn(\beta)}^2}
d\beta$;  then $\tau$ is invertible with inverse
$\tau\mapsto \alpha(\tau)$, $\tau \in I$.
We work with the parametrization
$$
\tau\mapsto x_*(\tau)= \widetilde x(\alpha(\tau))
$$
and  then
$$x_*'(\tau)=
(x_2(\alpha(\tau))),- x_1(\alpha(\tau)).
$$

In view of an analogue of  \eqref{div} it suffices to show that
\begin{equation}
\Big|\int_I e^{-\ic\inn{x_*(\tau)}{\xi}} \chi(\tau) d\tau
\Big|\lc (\log(2+|\xi|))^{-\ga}
\label{eq:5.7}
\end{equation}
Let $c_0= r_1/2r_2$ and $g(\tau, \xi)=\inn{\tfrac{\xi}{|\xi|}}{\tfrac {x_*'}{|x_*'|}}$.
Fix $|\xi|\ge 2$.
We split  our interval $I$ into no more than $16$ subintervals, where on each interval $J$
either $|g(\tau,\xi)|\ge c_0/4$ for all $\tau \in J$ or
 $|g(\tau,\xi)|\le c_0/4$ for all $\tau \in J$; and
$\tau\mapsto g(\tau,\xi)$ is  monotonic on $J$.

Suppose for all $\tau \in J$ we have $|g(\tau, \xi)|\ge c_0/4$. Then
$|\inn{x_*'(\tau)}{\xi}|\ge c_1|\xi|$ in $J$ and by van der Corput's Lemma we get
\begin{equation}
\Big|\int_J e^{-\ic\inn{x_*(\tau)}{\xi}} \chi(\tau) d\tau\Big|\lc |\xi|^{-1}
\label{eq:5.8}
\end{equation}
which of course is much better then the desired estimate.

Now fix $J'$ with the property that  $|g(\tau,\xi)|\le c_0/4$ for all $\alpha\in J'$.
Now $x(\alpha(\tau))$ and $x_*'(\alpha(\tau))$ are orthogonal and
thus $|\inn {\tfrac x{|x|}}{\tfrac \xi{|\xi|}}|\ge (1-c_0^2/16)^{1/2}\ge 1-c_0/4$; moreover
\begin{eqnarray*}
\big|\inn{\fn}{\tfrac{\xi}{|\xi|}}\big|
&\ge &
\big|\inn{\tfrac x{|x|}}{\fn}  \inn{\tfrac x{|x|}}{\tfrac\xi{|\xi|}}\big|-
\big|\inn{\tfrac {x_*'}{|x_*'|}}{\tfrac{\xi}{|\xi|}}\big| \\
&\ge & c_0(1-c_0/4)-c_0/4\ge c_0/4.
\end{eqnarray*}
Now for $\alpha \in J'$ we have $\fn(\alpha)=(x_2'(\a), -x_1'(\a))$ and thus
$$|\inn{ x_*''(\tau)}{\xi}|=|\inn{\fn(\alpha(\tau))}{\xi} \a(\tau)|\ge \frac{c_0}4 |\alpha'(\tau)|.
$$
Therefore
\begin{eqnarray*}
\int_{J'} \log^\ga(2+ \frac1{
|\inn{ x_*''(\tau)}{\xi/|\xi|}|})^{\ga} d\tau
&\ge & c_1\int_{J'} \log^\ga(2+ \frac1{|\alpha'(\tau)|})
 d\tau \\
& = & c_1\int_{\alpha(J')} \log^\ga(2+ |\tau'(\alpha)|) |\tau'(\alpha)| d\alpha
 d\tau
\end{eqnarray*}
where $c_1$ is chosen independently of $\xi$.
The latter expression is finite since $|\tau'(\alpha)|\approx |\kappa(\alpha)|$
which is assumed to be in $L\log^{\ga}\! L$.
Thus we may apply Lemma \ref{lem:5.2} with $\la=|\xi|$ and obtain
\begin{equation}
\Big|
\int_{J'} e^{-\ic\inn{x_*(\tau)}{\xi}} \chi(\tau) d\tau\Big|\lc
(\log(2+|\xi|))^{-\ga}
\label{eq:5.9}
\end{equation}
and the assertion follows from combining
the estimates \eqref{eq:5.8} and \eqref{eq:5.9} on the various subintervals.
\qed

\medskip

\noindent{\bf Proof of Proposition \ref{prop:5.1}}.
Given Lemma \ref{lem:5.3}  this is just an application of the standard
argument.
 Let $N_\eps^*(t)$, $E_\eps^*(t)$ be defined as in \eqref{Neps} and
 \eqref{Eeps}, but for the set $\Omega^*$ in place of $\Omega$.
Then by Lemma \ref{lem:5.3} for $t\ge 2$
\begin{eqnarray*}
|E_{\eps}^*(t)| & \le & (2\pi t^2)\sum_{k\neq 0}|
\widehat{\chi_{\Om^*}}(2\pi tk) |
|\widehat\zeta(2\pi \eps k)|\\
&\lc &  t^2\sum_{k\neq 0}(t|k|)^{-1}\log(2+|t|k|)^{\ga}(1+\eps|k|)^{-N}\lc
t (\log t)^{-\ga}\eps^{-1}.
\end{eqnarray*}
Also $N_\eps^*(t-C\eps)\le N^*(t)\le N_\eps^*(t+C\eps)$ and
applying the previous estimate with
$t\pm C\eps $ in place of $t$ yields
$$
|E_{\Om^*}(t)|\lc \big[ t (\log t)^{-\ga}\eps^{-1}+ t\eps\big].
$$

Thus for the choice $\eps=(\log t)^{-\ga/2}$ we obtain the
asserted estimate.\qed

\section{A weighted estimate for lattice points \\on lines in thin annuli}
\label{sectiondiscrete}

The purpose of this section is to prove  a bound for sums
$\sum_\ell \mu(\lunit,\tfrac{1}{R|\ell|})$ where the sum
runs over the lattice points contained on a given line segment
in the $\rho^*$-annulus
\begin{equation}\label{rhostarannulus}
\mathcal{A}\left( r,h\right) =\left\{ x\in \mathbb {R}^{2}:r\leq
\rho ^{*}\left(
x\right) \leq r+h^{-1}\right\} .
\end{equation}

It turns out that if $h\gc R^{1/2}$ and the line segment is
 sufficiently long then the trivial bound
$\mu(\lunit,\tfrac{1}{R|\ell|})\lc (Rr)^{-1/2}$ may be substantially
improved for most of the lattice points on the line segment; {\it i.e.}
the fact that the thin $\rho^*$ annulus contains long line
segments reflects a rather fast  turning of the normals for the
original domain
which may only happen if $\partial \Omega$ lacks smoothness.

Throughout this and the next two chapters we shall adopt the following
notations.
Given certain subsets $A$, $B$, $G$,  $I$, $J$ etc.
we shall use blackboard bold fonts to denote  by
$\bbA$, $\bbB$, $\bbG$, $\bbI$, $\bbJ$
the intersections of these sets with the integer lattice $\bbZ^2$
(the standard notation for the plane $\mathbb  R^2$ remains an exception
to this convention).
We shall
adopt  the convention that a \emph{line segment} $I=%
\overrightarrow{PQ}$ is a nontrivial segment whose endpoints $P,Q$ lie in
the lattice $\mathbb {Z}^{2}$.
The corresponding
collection $\bbI=I\cap \bbZ^{2}$ of lattice points in $I$ will be
called a \emph{lattice line segment}.

 If $n$ is an odd natural number and $J$ is a line
segment, we let $nJ$ denote the line segment concentric with
and parallel to
$J$ but with $n$ times the length.
The distance between
 consecutive lattice points on the line segment $J$  is constant. Let
$d\equiv d(J)$
denote this distance; then
$d^{-1}=(\card\, \bbJ-1)/|J|$
is the density of
lattice points on $J$.

The following result is the key for the
proof of  Proposition \ref{prop:3.3}. We are seeking to improve the
bound from \eqref{eq:1.7} (with $\delta= (R|\ell|)^{-1}$)
\begin{equation}\label{triviallinebound}
\sqrt{Rr}\sum_{\ell \in \mathbb {I}}\mu(\lunit,\tfrac{1}{R|\ell|})
\lc \card (\bbI),
\end{equation}
if $\card (\bbI)$ is sufficiently large.

\begin{lemma}\label{line}
Let $\Omega $ be an open convex bounded set in the plane $\bbR^2$,
with positive curvature, and containing the origin. Let $10\leq r\leq
R<\infty $, and let $h\ge R^{1/2}$.
Let $J$ be a closed line segment whose endpoints are contained in
$\bbZ^2$ and let  $\mathbb {J}=J\cap \mathbb {Z}^{2}$.
Assume that the ninefold dilate $9J$ is  contained in $\cA(r,h)$ and that
$\card(\bbJ) \ge 10 $.

Let $d=d(\bbJ)$ and \begin{equation} \label{criticalparameter}
\cT(\bbJ)\equiv \cT(R,\bbJ,h,r): = R^{3/4} d(\bbJ)^{1/2}h^{-1/2} r^{-1/4}.
\end{equation}
Then
\begin{equation}\label{claimoflemma61}
\sqrt{Rr}\sum_{\ell \in \mathbb {J}}\mu(\lunit,\tfrac{1}{R|\ell|})
\lc  \begin{cases}
 \Big(\frac{Rr}{h^2d^2}\Big)^{1/4} &\quad\text{ if }
\card(\bbJ)\le \cT(\bbJ),
\\
\Big(\frac{r\card(\bbJ)}{hd^2}\Big)^{1/3}
&\quad\text{ if }
\card(\bbJ)\ge \cT(\bbJ).
\end{cases}
\end{equation}


\end{lemma}

\noi{\bf Proof.}


We begin by introducing some additional notation.
For $P\in \mathbb {R}^{2}$, let $\fn_{P}^{*}$ denote the unit
outward normal to $\Omega ^{*}$ at the point $\frac{P}{\rho ^{*}\left( P\right) }
\in \partial \Omega ^{*}$. Note that this
unit normal is uniquely determined since it is parallel to the uniquely
determined position vector $\overrightarrow{O\mathcal{P}_{P}}$ joining the
origin $O$ to the point $\mathcal{P}_{P}\in \partial \Omega $ having
$P$ as
one of its outward normal vectors.



%

We set
 $I=3J$ and let
$P$, $Q$  be the endpoints of
$I$.
Our first observation is that the  angle
between the vectors
$\fn_{P}^{*}$ and $\fn_{Q}^{*}$ (the nonnegative angle less than
$\pi $)
satisfies
\begin{equation} \label{phi0}
\psi:=\measuredangle \left(\fn_{P}^{*},\fn_{Q}^{*}\right)
 \leq C_1\arctan\Big( \frac{h^{-1}}{| I|}\Big) ,
\end{equation}
where $| I| =| \overline {QP}| $ is the length of the
segment $I$ and $C_1$ is a constant depending only on $\Omega$.
To see this one notes that there is a rectangle of width  $ch^{-1}$
which contains the  line segment $3I$ and which is contained in the
annulus $\cA(r,h)$.  If we had
$\tan\measuredangle \left(\fn_{P}^{*},\fn_{Q}^{*}\right)
 \geq \widetilde Ch^{-1}| I|^{-1}$ for large $\widetilde C$ then
it is easy to see that the triangle $OPQ$ would contain points in the complement of $\{\rho^*\le r+h\}$. Thus we have \eqref{phi0}.

For $x \in I$, define the
collection $\Gamma( x) $ of subsegments
$J'=\overrightarrow{UV}$ of $I$ by
\[
\Gamma( x) =\big\{ J'=\overrightarrow{UV}\subset
I:x \in J', |J'| \geq d(\bbJ)\big\} ,
\]
and the corresponding uncentered maximal function $\mathcal{M}$ on
$I $ by
\[
\mathcal{M}\left( x\right) =\sup_{J'=\overrightarrow{UV}\in
 \Gamma(x) }\frac{1}{|J'| }\measuredangle \big(
\fn_{U}^{*},\fn_{V}^{*}\big) .
\]
Then
$\mathcal{M}$ is in  weak $L^1$ on $I$
by F. Riesz's lemma (\cite{riesznagy}, ch. 1); it
 satisfies the inequality
\begin{equation}\label{weak11}
\left| \left\{  x \in I:\mathcal{M}( x)
>\lambda \right\} \right| \leq 2\psi/\lambda .
\end{equation}

We shall now consider a decomposition of the set $\bbI$
 depending on a parameter $q$; we shall see that the choice
\begin{equation}
\label{qchoice}
q= \begin{cases}\Big(\frac{Rrd^2}{h^2}\Big)^{1/4} &\text{ if }
\card(\bbJ)\le \cT(\bbJ)
\\
\Big(\frac{rd\,\card(\bbJ)}{h}\Big)^{1/3}
&\text{ if }
\card(\bbJ)\ge \cT(\bbJ)
\end{cases}
\end{equation}
will be (essentially) optimal.

Let $\mathbb {B}\subset \bbI:=I\cap\bbZ^2$ denote the set
\[
\mathbb {B}=\big\{ \ell \in \mathbb {I}:\mathcal{M}
( \ell)
> q^{-1}\psi \big\}
\]
to which, following the terminology in Calder\'on-Zygmund theory,
we refer  as the
set
of \emph{bad} lattice points. Let $\mathbb {G}=\mathbb {I}
\setminus \mathbb {B}$ be the set of \emph{good}
lattice points. Denote by $u$ a unit vector parallel to $I$.
 Then the segments $\overrightarrow{(\ell-du/2)(\ell+du/2)}$
 are pairwise disjoint and for each $\ell \in \mathbb {B}$, either
$\overrightarrow{( \ell -du/2)\ell}$ or
$\overrightarrow{\ell ( \ell +du/2) }$ is
contained in $$\big\{ x \in I:\mathcal{M}(x) >
\psi /q\big\} .$$ Thus by (\ref{weak11}) we have
\begin{align*}
\card(\mathbb {B}) &= d^{-1}\Big| \bigcup _{\ell \in \mathbb {B}}
\overrightarrow{( \ell -du/2) ( \ell +du/2)}\Big| \\
&\leq 2d^{-1}\big|\big\{ x \in I:\mathcal{M}( x) >
\psi/q \big\} \big| \lc q/d.
\end{align*}
By \eqref{eq:1.7} we obtain
\begin{equation}
\label{badestimate}
\sum_{\ell\in \bbB}\sqrt{Rr}\mu(\lunit,\tfrac{1}{R|\ell|})
\lc q/d
\end{equation}

We shall  now obtain an estimate for the  sum over $\ell \in
\mathbb {G}\cap \mathbb {J}$
(rather than over \emph{all} of $\mathbb {G}$).
Note that if $\ell \in \mathbb {G}\cap \mathbb {J}$, then
\begin{eqnarray}
\measuredangle \left( {\fn}_{\ell -\alpha du}^{*},
{\fn}_{\ell +\alpha du}^{*}\right) &\leq &2\alpha d \psi/q, \\
\measuredangle \left(  \ell -\alpha du,\ell\right ),\,
\measuredangle \left(\ell , \ell +\alpha du\right)
&\geq &c\alpha d/r
\end{eqnarray}
for $1/2\le \alpha \leq | J|/d$ upon using
$| \ell| \approx r$.

 Passing to the dual set $\Omega^{**}=\Omega $ with
defining Minkowski  functional $\rho $, we have that for $\alpha \geq 0$,
the points
$\mathcal{A}_{\ell, \alpha }^\pm=\frac{\fn_{\ell \pm\alpha du}^{*}}{\rho
( \fn_{\ell \pm\alpha du}^{*}) }$  in $\partial \Omega $
have unit normals  $\mathbf{\nu }_{\ell,\alpha }^\pm=\frac{
\ell \pm\alpha du}{\left| \ell \pm\alpha du\right| }$,
 and that
\begin{equation}\label{estimates}
\begin{aligned}
\big|\measuredangle
 \big( {\nu }_{\ell,\alpha }^\pm,\lunit\big)\big| &\geq c\alpha d/r
\\
|\measuredangle \left( \mathcal{A}_{\ell,\alpha }^-,
\mathcal{A}_{\ell,\alpha }^+\right)|
&\leq 2\alpha d \psi/q
\end{aligned}
\end{equation}
for $1/2\leq \alpha \leq | J|/d$.

Using this bound we derive a bound on the diameter
of the cap
$\mathcal{C}_{\Omega }( \tfrac{\ell }{| \ell | },
\tfrac{1}{R| \ell| }) $. Let
\[
D=2C_1d  (h|I|q  )^{-1},
\]
where $C_1$ is as in \eqref{phi0}.
 Suppose that $\partial \Omega $ is parametrized in a neighborhood of
$\frac{\fn_{\ell }^{*}}{\rho \left( \fn_{\ell }^{*}\right) }$
by
$t\mapsto \gamma(t)$, $|t|\le c'$,
with
$$\inn{\gamma(t)-\cA_{\ell,0}}{\lunit}= \varphi_\ell(t);$$
here
$\cA_{\ell,0}^\pm
=\frac{\fn_{\ell }^{*}}{\rho \left( \fn_{\ell }^{*}\right) }$,
 and $\varphi_\ell(t)$ is convex and nonnegative with $\varphi_\ell(0)=0$,
$\varphi_\ell'(0)=0$.

Now by (\ref{phi0}) we have $D\geq 2d \psi/q$,
and so (\ref{estimates}) shows that
\[
| \varphi_\ell ^{\prime }( \alpha D) |
\geq
| \varphi_\ell ^{\prime }( \alpha 2d\psi/q ) |
 \geq c\alpha
d/r,
\;\;\;\;\;
1/2\leq \alpha \leq |J|/d,
\]
and it follows that for $D/2\leq |t|\leq|J|d^{-1}  D$,
\begin{equation}\label{boundbelowvarphi}
\varphi_\ell \left( t\right) \geq \int_{D/2}^{t}\frac{cd}{Dr}s\,ds=
\frac{1}{2}\frac{cd}{Dr}\left[ t^{2}-\left( D/2\right)
^{2}\right] .
\end{equation}
Now assume $|J|/d\ge 1$ and observe that then
for $t=D|J|/d$, the upper bound for
 $|t|$, we have by \eqref{boundbelowvarphi} and the definition of $D$
and $J$
\begin{align*}
&\varphi_\ell \big( \tfrac{| J| }{d}D\big)
\ge\frac{c}{2}\frac{dD}{r}\big[(|J|/d)^2-1/4\big]\geq \frac{3c|J|^2 D}{8dr}
= \frac{cC_1 |J|}{4rhq}.
\end{align*}
Since $h\le R$
and in view of $h\ge R^{1/2}$, $1\le r\le R$ we see that
the choice of $q$ in \eqref{qchoice} certainly implies
$$
\varphi_\ell ( \tfrac{\left| J\right| }{d}D)\ge c'(Rr)^{-1}\ge c''
(R|\ell|)^{-1}
$$
This estimate  and \eqref{boundbelowvarphi} imply that
\begin{align*}
\mu(\lunit, \tfrac{1}{R|\ell|}) &\le C(c'')
\mu(\lunit, \tfrac{1}{c''R|\ell|})
\lc\varphi_\ell ^{-1}(
\tfrac{1}{c''R| \ell | })
\\&\lc\varphi_\ell ^{-1}(
\tfrac{1}{R| \ell | })
 \leq D+2\sqrt{\frac{2Dr}{cd}
\frac{1}{R\left| \ell \right| }} \\
&\lc ( D+(D/dR)^{1/2}),
\end{align*}
and since
$D=2C_1 d   (h| I|q)^{-1}$
we then have
$$
\sqrt{Rr}\mu(\lunit,\tfrac{1}{R|\ell|})\lc
\frac{R^{1/2}}{ h} \frac{dr^{1/2}}{|I|q }
+
\Big(\frac{r}{h|I|q}\Big)^{1/2},
$$
for $\ell\in \bbG\cap\bbJ$.
Since $|I|\approx d\,\card(\bbI)\approx d\, \card (\bbJ)$ we obtain
\begin{equation}\label{goodestimate}
\sum_{\ell\in \bbG\cap \bbJ}
\sqrt{Rr}\mu(\lunit,\tfrac{1}{R|\ell|})\lc \frac{(Rr)^{1/2}}{hq}+
\Big(\frac{r\card(\bbJ)}{hdq}\Big)^{1/2}.
\end{equation}
Altogether
\eqref{goodestimate} and \eqref{badestimate} yield
\begin{equation}\label{goodbadestimate}
\sum_{\ell\in \bbJ}
\sqrt{Rr}\mu(\lunit,\tfrac{1}{R|\ell|})\lc \frac{q}{d}+
 \frac{(Rr)^{1/2}}{hq}+
\Big(\frac{r\card(\bbJ)}{hdq}\Big)^{1/2}.
\end{equation}

We essentially choose $q$ to minimize this expression.
Its square  is comparable to $d^{-2} F(q)$ where
\begin{equation}\label{Fq} F(q)= q^2+b_1 q^{-2} + b_2 q^{-1},
\end{equation}
and the positive coefficients are given by
$$b_1= Rr d^2 h^{-2}, \qquad b_2=
r\card(\bbJ) dh^{-1} .$$ In what follows we shall need the relation
$$b_2^{4/3}/b_1= \big( \card(\bbJ)/\cT\big)^{4/3}$$
where $\cT$ is as in  \eqref{criticalparameter}.
To analyze \eqref{Fq}  it is natural to distinguish two cases,
where in the first case the
second term $b_1 q^{-2}$
dominates the   third term $b_2 q^{-1}$
and in the second
case we have the opposite inequality.

In the first case, $q\le b_1/b_2$, we have that
$F(q)\approx q^2+b_1q^{-2}$ and the latter expression has a local
minimum where $q=b_1^{1/4}$. This value $b_1^{1/4}$ belongs
to the currently relevant interval $(0,b_1/b_2]$
 when
$b_2\le b_1^{3/4}$ which is equivalent to $\card(\bbJ)\le \cT$.
In this case we thus make the choice $q=b_1^{1/4}$
(which is the value $(Rrh^{-2}d^{-2})^{1/4}$  in \eqref{qchoice}).
Then the right hand side of  \eqref{goodbadestimate}
is bounded by a constant times
\begin{equation}
d^{-1}\sqrt{ F(b_1^{1/4})} \lc
d^{-1} b_1^{1/4}=(Rr d^{-2} h^{-2})^{1/4}.
\end{equation}

In the second case, $q\ge b_1/b_2$ we have
  $F(q)\approx q^2+2b_2 q^{-1}$ and the latter expression has a local
minimum at $b_2^{1/3}$ which lies in the interval
$[b_1/b_2,\infty)$  when $b_2^{4/3}\ge b_1$; this condition is
equivalent with $\card(\bbJ)\ge \cT$.
Thus we now make the choice
$q=b_2^{1/3}=(r\card(\bbJ) d h^{-1})^{1/3}$.
Now the right hand side of  \eqref{goodbadestimate}
is bounded by a constant times
\begin{equation}
d^{-1}\sqrt{ F(b_2^{1/3})} \lc
d^{-1} b_2^{1/3}= (r\card (\bbJ) h^{-1} d^{-2})^{1/3}
\end{equation}
and the
estimate  \eqref{claimoflemma61} is established.
\qed

\section{An  elementary lemma}

We give a standard estimate
on lattice points on convex polygons.

\begin{lemma}\label{twothirds}
Suppose we are given integer points $F_1,\dots, F_J$ in $\bbZ^2$,
which are vertices of a convex polygonal curve; i.e.
the interiors of the line segments
$\overline {F_iF_{i+1}} $ are mutually disjoint and if
$\overrightarrow{F_{i}F_{i+1}}=L_i (\cos \beta_i, \sin \beta_i)$, $i=1,\dots, J-1$,
then we have $L_i>0$ and $\beta_{i+1}>\beta_i$. Suppose also that
$\beta_{J-1}-\beta_1\le 2\pi$, and set $L=\sum_{i=1}^{J-1}L_i$.
Then $$J\le 2+(\beta_{J-1}-\beta_1)^{1/3}L^{{2}/{3}}.$$
\end{lemma}

\noi{\bf Proof.}
Let $\Delta(A,B,C)$
 denote the triangle with vertices $A,B,C$. Then we use again that the
area
of $\Delta(A,B,C)$ is at least $1/2$ if $A,B,C$ are noncollinear lattice points.
Thus
\begin{align*}
J-2 &\leq \sum_{j=1}^{J-2}\big(2\;\area( \Delta (
F_{j},F_{j+1},F_{j+2}) ) \big)^{1/3}
\\
&=
\sum_{j=1}^{J-2}\big(| \overrightarrow{F_{j}F_{j+1}}|\,
|\overrightarrow{F_{j+1}F_{j+2}}|\,|\sin (\beta_{j+1}-\beta_j)|\big)^{1/3}
\end{align*}
and thus by H\"older's inequality
\begin{align*}
J-2&\leq \Big(\sum_{j=1}^{J-2}| \overrightarrow{F_{j}F_{j+1}}|
\Big)^{1/3}\Big( \sum_{j=1}^{J-2}
| \overrightarrow{F_{j+1}F_{j+2}}|\Big)^{1/3}
\Big(\sum_{j=1}^{J-2}|
\sin(\beta_{i+1}-\beta_i)| \Big)^{1/3}\\
&\leq L^{2/3}(\beta_{J-1}-\beta_1)^{1/3}. \qed
\end{align*}

\medskip

In particular if $\fP$ is a convex polygon of length $L$ whose
vertices are integer lattice points then the number of vertices of $\fP$
 is $O(L^{2/3})$.
This is a special case of Andrews' result \cite{Andrews}.

\section{
Proof of Proposition \ref{prop:3.3}}\label{proofprop33}

In what follows we fix $r,R,h$ with $10\le r\le R$, $R^{1/2}\le h\le R$.
Let
\[
\Omega _{\pm }^{*}=\left\{  x \in \Bbb{R}^{2}:\rho ^{*}(x)
<r \pm h^{-1}\right\}
\]
and $\mathcal{A}=\Omega _{+}^{*}\setminus \Omega _{-}^{*}$.
Denote by $\Bbb{A}=\mathcal{A}\cap \Bbb{Z}^{2}$
the set of lattice points in the annulus $\mathcal{A}$.
Let $\mathcal{H}$ be the convex hull of $\Bbb{A}$ and let
 $\Bbb{E}=\left\{ E_{j}\right\} _{j=1}^{\card(\bbE)}$
be the extreme points of $\mathcal{H}$ arranged in counterclockwise order
around the origin. First observe
that from Lemma
\ref{twothirds} we have
\begin{equation} \label{initialboundforJ}
\card (\Bbb{E})\lc r^{2/3}.
\end{equation}

Recall our convention from the
\S\ref{sectiondiscrete} that a \emph{line segment}
is a nontrivial segment $I$ whose endpoints lie in the
lattice $\Bbb{Z}^{2}$.
Our second observation is that every lattice point in $\bbA\setminus
\bbE$ belongs to
some line segment $I$ that contains an extreme point from $\Bbb{E}$ and lies
entirely within the annulus $\mathcal{A}$. More precisely, let $\ell \in
\Bbb{A}$ and suppose that for some $1\leq j\leq \card(\bbE)$, $\ell $ belongs to the
closed triangular sector $\mathcal{S}_{j}=\Delta \left(
E_{j},E_{j+1},0\right) $ with vertices $E_{j}$, $E_{j+1}$ and the origin.
Then the convex set $\Omega _{-}^{*}$ cannot intersect both of the line
segments $\overrightarrow{E_{j}\ell }$ and $\overrightarrow{E_{j+1}\ell }$,
and hence at least one of them must lie in $\mathcal{A}$.

Thus for $1\leq j\leq \card(\bbE)$, $\Bbb{A}_{j}=\Bbb{A}\cap \mathcal{S}_{j}$ is
contained in the union $\mathcal{I}_{j}$ of all maximal line segments $I$
contained in $\mathcal{A}\cap \mathcal{S}_{j}$ having one endpoint that is
either $E_{j}$ or $E_{j+1}$. We further distinguish the segments $I$ in
$\mathcal{I}_{j}$ by letting $\{ I_{j,n}^{-}\}$, $n=1,\dots N_{j}^{-}$
be an enumeration of the line segments in $\mathcal{I}_{j}$ having $E_{j}$
as an endpoint, and arranged clockwise about $E_{j}$ as $n$ increases from $
1 $ to $N_{j}^{-}$. Similarly, $\{ I_{j,n}^{+}\}$, $n=1,\dots,N_{j}^{+}$
 is an enumeration of the line segments in $\mathcal{I}_{j}$
having $E_{j+1}$ as an endpoint, and arranged counterclockwise about $
E_{j+1}$ as $n$ increases from $1$ to $N_{j}^{+}$.
Note that if
the line segment $\overrightarrow{E_{j}E_{j+1}}$
joining the consecutive extreme points $E_{j}$ and
$E_{j+1}$ is contained in $\cA$ then
$I_{j,N_{j}^{+}}^{+}=I_{j,N_{j}^{-}}^{-}=\overrightarrow{E_{j}E_{j+1}}$.

The next lemma implies that the total number of line segments in
$$\mathcal{I}:
=\bigcup _{j=1}^{\card(\bbE)}\mathcal{I}_{j}$$ still does not exceed $C r^{2/3}$
 if we assume $r\le R$ and $h\ge R^{1/2}$.
For this we define $L_{j}=| \overrightarrow{E_{j}E_{j+1}}| $ to
be the length of the segment $\overrightarrow{E_{j}E_{j+1}}$, and $\Theta
_{j}$ to be the (positive) change of angle for the
 normal direction to $\partial \Omega^*$ across the sector $\mathcal{S}_{j}$
(recall that every point in $\partial \Omega^*$ has a unique  normal since
$\Omega$ has positive curvature).
Specifically

\begin{lemma}
\label{segmentbound}For each $1\leq j\leq \card(\bbE)$, we have
\[
N_{j}^{-}+N_{j}^{+}\leq C\big( 1+\area( \mathcal{A}\cap \mathcal{S}%
_{j}) +L_{j}^{2/3}\Theta _{j}^{1/3}\big).
\]
\end{lemma}

\proof  We will establish the indicated estimate for $%
N_{j}^{+}$, the case for $N_{j}^{-}$ being similar. Fix $j$ and let the
segment $I_{j,n}^{+}$ have endpoints $E_{j+1}$ and $F_{n}$ so that $%
I_{j,n}^{+}=\overrightarrow{E_{j+1}F_{n}}$. Consider the collection of
segments $\{ \overrightarrow{F_{n}F_{n+1}}\}
_{n=1}^{N_{j}^{+}-1}$ and set
\begin{eqnarray*}
\fT_j &=&\big\{ n:1\leq n<N_{j}^{+}\text{ and }\overrightarrow{F_{n}F_{n+1}}
\subseteq \mathcal{A}\cap \mathcal{S}_{j}\big\} , \\
\fP_j &=&\big\{ n:1\leq n<N_{j}^{+}\text{ and }\overrightarrow{F_{n}F_{n+1}}
\nsubseteq \mathcal{A}\cap \mathcal{S}_{j}\big\} .
\end{eqnarray*}

We first note that  the triangles $\Delta(E_{j+1},F_{n},F_{n+1}) $
are pairwise disjoint and contained in
$\mathcal{A}\cap \mathcal{S}_{j}$ for $n\in \fT_j$.
As the corners of these triangles are in $\bbZ^2$ the area of each
triangle is at least $1/2$ and we estimate
\begin{equation}\label{triangularsegments}
\card(\fT_j) \leq 2\sum_{n\in \fT_j}\area\big( \Delta (
E_{j+1},F_{n},F_{n+1} )\big)
\leq 2\;\area( \mathcal{A}\cap \mathcal{S}_{j}) .
\end{equation}
Now the integers in $\fP_j$ can be written as a union of maximal chains $
\mathcal{C}_{i}$ of consecutive integers as follows:
\[
\fP_j=\cup _{i=1}^{M_j}\mathcal{C}_{i},\;\;\;\;\;\mathcal{C}_{i}=\left\{ n\right\}
_{n=a_{i}}^{b_{i}},
\]
where $a_{i}\leq b_{i}$ and $b_{i}+2\leq a_{i+1}$. Note that the number of
chains $M_j$ satisfies
\[
M_j\leq 1+\card(\fT_j)\leq 1+
2\;\area\left( \mathcal{A}\cap \mathcal{S}_{j}\right) .
\]

For each $n\in \fP_j$ we may write
$\overrightarrow{F_{n}F_{n+1}}
=\rho_n(\cos\alpha_n, \sin \alpha_n)$ where $\rho_n>0$ and
$\alpha_n>\alpha_m$ if $m>n$. In particular the lines associated to a
fixed chain form a convex polygon.
To see this, let $m>n$ and  note that the convex set $\Omega_-^*$
intersects
both of the
line segments $\overline{F_nF_{n+1}}$ and $\overline{F_mF_{m+1}}$,
and that the intersection with $\overline{F_mF_{m+1}}$
 occurs to the clockwise side of the intersection with
$\overline{F_nF_{n+1}}$ (we adopt an obvious convention
regarding the use of the phrase "to the clockwise side of").
Thus the angle
$\alpha_m$ of $\overline{F_mF_{m+1}}$ is less than the angle $\alpha_n$ of
$\overline{ F_nF_{n+1}}$.
Moreover the sectors generated by $\{O,F_n, F_{n+1}\}$
and $\{O,F_m, F_{m+1}\}$ have disjoint interior so that
\[\sum_{n\in \fP_j}\big|\overrightarrow{F_{n}F_{n+1}}\big| \lc
L_j.\] See the following diagram for a somewhat artificial 
illustration with $\frak P_j = \{1,2\}$ and $\frak T_j = \{3\}.$
Note that the segment  $\overline{F_2F_3}$ lies to the
 clockwise side of the segment $\overline{F_1F_2}$, and that $F_4 = E_j$. 
Realistic illustrations with actual lattice points require
extremely thin annuli  and are difficult to  draw.
\begin{center}\label{diagram}
\includegraphics[
height=3in, width=5in]
 {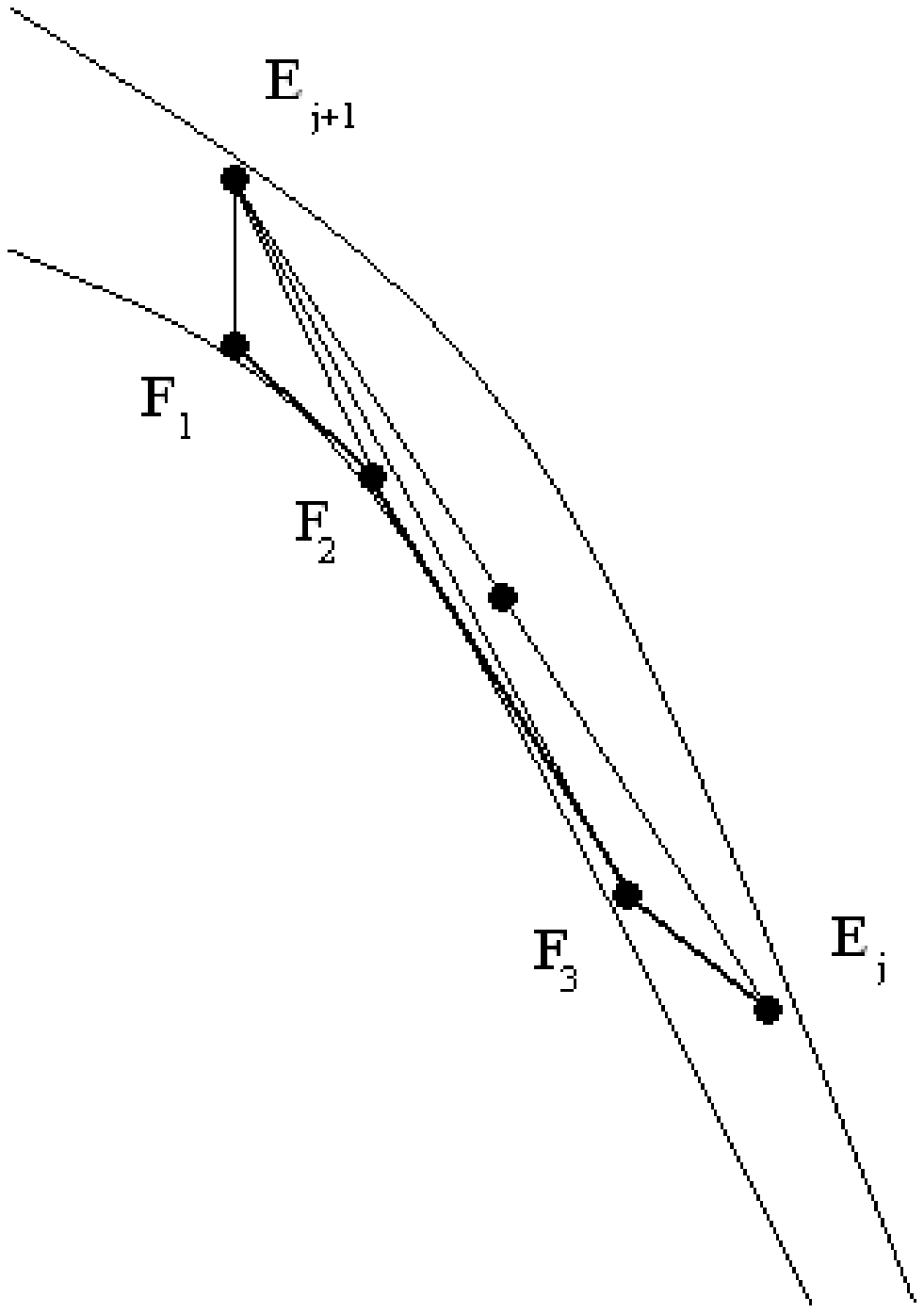}
\end{center}
Now consider a chain $\mathcal{C}_{i}=\{ n\} _{n=a_{i}}^{b_{i}}$
of length $b_{i}-a_{i}+1\geq 3$ and the associated set of lattice
points $\{ F_{n}\} _{n=a_{i}}^{b_{i}}$.
By Lemma \ref{twothirds}, we have
\[
b_{i}-a_{i}+1\leq 2+\Big( \sum_{n=a_{i}}^{b_{i}-1}\left| \overrightarrow{%
F_{n}F_{n+1}}\right| \Big) ^{\frac{2}{3}}\big(
\sum_{n=a_i}^{b_i-1} (\alpha_{n+1}-\alpha_n)
\big)^{\frac{1}{3}}
\]
for all chains of length $b_{i}-a_{i}+1\geq 3$, and also trivially for
chains of length $1$ or $2$ as well. Summing in $i$ from $1$ to $M_j$ we thus
obtain
\begin{align*}
\card(\fP_j)
&=\sum_{i=1}^{M_j}( b_{i}-a_{i}+1) \\
&\leq 2M_j+\sum_{i=1}^{M_j}\Big\{ \big( \sum_{n=a_{i}}^{b_{i}-1}|
\overrightarrow{F_{n}F_{n+1}}| \big)^{\frac{2}{3}}
\big(
\sum_{n=a_i}^{b_i-1} (\alpha_{n+1}-\alpha_n) \big)^{\frac{1}{3}} \Big\}
\end{align*}
and thus by H\"older's inequality
\begin{align*}
\card(\fP_j)&\le 2M_j+\Big( \sum_{i=1}^{M_j}\sum_{n=a_{i}}^{b_{i}-1}\big|
\overrightarrow{F_{n}F_{n+1}}\big| \Big) ^{\frac{2}{3}}\Big(
\sum_{i=1}^{M_j}
\sum_{n=a_i}^{b_i-1} \alpha_{n+1}-\alpha_n
\Big) ^{\frac{1}{3}} \\
&\le C\big(M_j +L_j^{2/3} \Theta_j^{1/3}
\big)\le C' \big(1+ \,\area(\cA\cap\cS_j)+
L_j^{2/3} \Theta_j^{1/3}\big).
\end{align*}
The Lemma follows by combining the inequalities for the cardinalities of $\fT_j$ and $\fP_j$.\qed

\medskip

We now proceed with the proof of Proposition
 \ref{prop:3.3}. First, by Lemma
\ref{segmentbound}
and inequality
\eqref{initialboundforJ} we can estimate the cardinality of $\cI$
by

\begin{align}
\card(\mathcal{I}) &\leq \sum_{j=1}^{\card(\bbE)}\card(\mathcal{I}_{j})\lc
\sum_{j=1}^{\card(\bbE)}\{ 1+\area( \mathcal{A}\cap \mathcal{S}_{j})
+L_{j}^{\frac{2}{3}}\Theta _{j}^{\frac{1}{3}}\}  \notag
\\
&\lc \Big( \card(\bbE)+\area ( \mathcal{A}) +\big(
\sum_{j=1}^{\card(\bbE)}L_{j}\big)^{\frac{2}{3}}\big(
\sum_{j=1}^{\card(\bbE)}\Theta
_{j}\big) ^{\frac{1}{3}}\Big) \notag
 \\
&\lc (r^{2/3}+r/h) \lc r^{2/3}
\label{numberofallsegments}
\end{align}
since $h\ge R^{1/2}\ge r^{1/2}$.


Now consider the lattice line segments  $\Bbb{I}_{j,n}^{\pm
}=I_{j,n}^{\pm }\cap \Bbb{Z}^{2}$  consisting of
the lattice points in $I_{j,n}^{\pm }$, and let $\Bbb A_j$ be the collection of lattice points which lie in $\cA\cap \cS_j$. We  then get
for each $j$
\begin{equation}
\Bbb{A}_{j}=\Bbb{A\cap }\mathcal{S}_{j}=\Big( \bigcup _{n=1}^{N_{j}^{+}}
\Bbb{I}_{j,n}^{+}\Big) \cup \Big( \bigcup _{n=1}^{N_{j}^{-}}
\Bbb{I}_{j,n}^{-}\Big)   \label{lattice}
\end{equation}
where $\sum_{j=1}^{\card(\bbE)}( N_{j}^{+}+N_{j}^{-}) \lc \card(\cI)\lc r^{2/3}$.

We now wish to apply Lemma \ref{line}, to the  intervals in $\cI$;
however
the assumption that the ninefold dilates are still contained in the annulus $\cA$ may not be satisfied. Therefore for every $I\in \cI$ we decompose $I$ in
subintervals $$I=\imath_+(I)\cup \imath_-(I)
\cup \bigcup_{m=-N(I)}^{N(I)}\imath_m (I)$$
where $9\imath_m(I)\subset I$.
Moreover if
$\imath_m(\bbI):=\imath_m(I)\cap \bbZ^2$,  and
$\imath_\pm(\bbI):=\imath_\pm(I)\cap \bbZ^2$,  and
then $\card (\imath_m(\bbI))\lc 2^{-|m|} \card(\bbI)$,
$\card(\imath_\pm(\bbI)) =O(1)$ and $N(I)\le C+\log_2(\card(\bbI))$.

We first have the trivial inequality
\begin{equation}\label{plusminustrivialbound}
\sum_{I\in \cI}\sum_\pm
\sum_{\ell \in \imath_\pm(\bbI)}
\sqrt{Rr}
\mu(\lunit,\tfrac{1}{R|\ell|}) \lc
 \card(\cI)\lc r^{2/3}.
\end{equation}

Now let $\fL$ denote the set of all lattice line segments $\{\imath_m(\bbI): |m|\le N(I), I\in \cI\}$.
We split $\fL$ into three subfamilies (here $\cT$ is as defined in
\eqref{criticalparameter}, i.e. $\cT(J)= R^{3/4} d(J)^{1/2} h^{-1/2} r^{-1/4}$).


\medskip
(i) $\fL_1$ consists of all
$\bbJ\in \fL$ which satisfy  $\card(\bbJ)\le \cT(J)$

\medskip
(ii) $\fL_2$ consists of all
$\bbJ\in \fL$ of the form $\imath_m(\bbI)$ for suitable $m$, $\bbI$,
where
$\card (\bbJ)>\cT(\bbJ)$ and $\card (\bbI)\le r^{1/3}$.

\medskip

(iii)
 $\fL_3$ consists of all
$\bbJ\in \fL$ of the form $\imath_m(\bbI)$ for suitable $m$, $\bbI$,
 where
$\card (\bbJ)>\cT(\bbJ)$ and $\card (\bbI)> r^{1/3}$.

\medskip

In our applications of Lemma  \ref{line}  we shall ignore the possible
nontrivial size of $d(J)$ and use the trivial bound $d(J)\ge 1$.

Notice that for all $I$ we have $N(I)\lc \log r$ so that
$\card (\fL)\lc r^{2/3}\log r$. Thus from the first inequality
in \eqref{claimoflemma61}  of
Lemma \ref{line} (and $d\ge 1$)
we get
\begin{equation}\label{Loneestimate}
\sum_{\bbJ\in \fL_1}\sqrt{Rr}\sum_{\ell \in \mathbb {\bbJ}}
\mu(\lunit,\tfrac{1}{R|\ell|}) \lc r^{2/3}\log r (Rrh^{-2})^{1/4}\lc
r^{11/12}\log r
\end{equation}
since $h^2\ge R\ge r$ (from the hypothesis in Proposition \ref{prop:3.3}).

Next we consider the lattice line segments in $\fL_2$  and use the second inequality  in \eqref{claimoflemma61}  of
Lemma \ref{line}. We obtain 
\begin{align}\label{Ltwoestimate}
&\sum_{\bbJ\in \fL_2}\sqrt{Rr}\sum_{\ell \in \mathbb {\bbJ}}
\mu(\lunit,\tfrac{1}{R|\ell|}) \notag\\
&\qquad\lc
\sum_{I\in \cI:\atop\card(\bbI)\le r^{1/3}}\sum_{m:\atop\imath_m(\bbI)\in \fL_2}
(rh^{-1}\, \card(\imath_m(\bbI))^{1/3}\notag
\\
&\qquad\lc
\sum_{I\in \cI:\atop\card(\bbI)\le r^{1/3}}
\big(rh^{-1}\, \card(\bbI)\big)^{1/3}\notag
\\&\qquad\lc  r^{1/3+1/9} \card(\cI) h^{-1/3
} \lc r^{17/18}
\end{align}
since $\card(\cI)\lc r^{2/3}$ and $h\ge R^{1/2}\ge r^{1/2}$.

Finally for the lattice line segments in $\fL_3$ we use again
 the second inequality  in \eqref{claimoflemma61}  of
Lemma \ref{line} but estimate differently

\begin{align}\label{Lthreeestimate}
&\sum_{\bbJ\in \fL_3}\sqrt{Rr}\sum_{\ell \in \mathbb {\bbJ}}
\mu(\lunit,\tfrac{1}{R|\ell|}) \notag\\
&\qquad
\lc
\sum_{I\in \cI:\atop \card(\bbI)> r^{1/3}}
\big(rh^{-1}\, \card(\bbI)\big)^{1/3}\notag
\\&\qquad\lc
\sum_{I\in \cI}
(r/h)^{1/3} \frac{\card(\bbI)}{r^{2/9}}\notag
\\
&\qquad\lc r^{1/3-2/9}h^{-1/3} \card (\bbA) \lc r^{17/18}
\end{align}
since $h\ge r^{1/2}$ and $\card (\Bbb A)\lc r$.

We combine estimates \eqref{plusminustrivialbound},
\eqref{Loneestimate},
\eqref{Ltwoestimate}, and
\eqref{Lthreeestimate} and deduce the asserted bound.\qed

\section{On the Sharpness of Proposition \ref{prop:3.3}}\label{sharpness}

Proposition \ref{prop:3.3} implies
the estimate
\begin{equation} \cK(R,h)\lc R^{-1}, \quad h\ge R^{1/2}
\label{ineqforKRH}
\end{equation} for the quantity defined in \eqref{bilinear}. We show that
the condition $h\ge R^{1/2}$ is needed in \eqref{ineqforKRH}.

More specifically we show that  there are positive
constants $c$ and $C$ such that for every $\varepsilon >0$ and $R>C$, there
exists an open convex bounded set $\Omega $ with curvature bounded uniformly below
(in the sense of any of the equivalent definitions in the subsequent section)
such that with $h=h(R)=
 R^{1/2}\left( \log R\right) ^{-\varepsilon }$ and suitable $C_0$
the quantity
$\cK(R,h)$ is
at least $c\varepsilon R^{-1} \log \log R$.

In what follows we let $r$ be a large  integer  satisfying
$$r\lc  (\log R)^{\epsilon/3}.$$
We will construct an open convex bounded set $\Omega $,
with curvature  $\gc 1/2$ everywhere
so that
\begin{equation}
\label{one}
(R|(m,n)|)^{1/2}
\mu( \tfrac{( m,n) }{| (m,n) | },\tfrac{1}{R| ( m,n) | }) =1
\end{equation}
for all  $( m,n) \in \Bbb{Z}^{2}$ with $0< |n| \leq m\leq r$,
and also so that
\begin{equation}
 \label{slim}
\left| \rho ^{*}\left( \left( m,n_{1}\right) \right) -\rho ^{*}\left( \left(
m,n_{2}\right) \right) \right|
 \le 2R^{-1/2}( \log R) ^{\varepsilon }
\end{equation}
for $0\leq \left| n_{1}\right| ,\left| n_{2}\right| \leq m\leq r$. With this
achieved, we restrict $k$ and $\ell $ in the sum on the left side of
\eqref{bilinear} to lie in the triangle of lattice points
$\Bbb{T}_{r}=\left\{ \left( m,n\right) :0\leq \left| n\right| \leq m\leq
r\right\} $. Writing $k=\left( m,n_{1}\right) $ and $\ell =\left(
m,n_{2}\right) $, we obtain
that with $h(R)=C_0^{-1}R^{\frac{1}{2}}\left( \log R\right)^{-\varepsilon }$,
\begin{align*}
\sum_{k,\ell \in \Bbb{Z}^{2}:| k| \leq R,| \ell |
\leq R,\atop
| \rho ^{*}( \ell ) -\rho ^{*}( k) |
\leq h(R)^{-1}}
&| k| ^{-2}\mu\big( \tfrac{k}{|k| },\tfrac{1}{R| k| }\big) \mu\big(
\tfrac{\ell }{| \ell | },\tfrac{1}{R| \ell | }\big)
\\&\geq R^{-1}\sum_{k,\ell \in \Bbb{T}_{r}:\atop
| \rho ^{*}( \ell )
-\rho ^{*}( k) | \leq h(R)^{-1}}| k| ^{-3}
\geq R^{-1}\sum_{(m,n)\in \bbT_r} m^{-2}\\&\ge c  R^{-1} \log r
= c R^{-1}\frac{\varepsilon }{3}\log \log R,
\end{align*}
the desired conclusion; here we have used  \eqref{one} and \eqref{slim},
respectively, in the first two inequalities above.

Now we give the details of the construction. Given $\varepsilon >0$, $r\in
\Bbb{N}$ large and $h>1$, denote by $\mathcal{R}_{m,n}$ the ray from the
origin $\left( 0,0\right) $ through the point $\left( m,n\right) $, and by
$K$ the circle of radius $hr^{2}$ centered at $\left(
1-hr^{2},0\right) $, so that  $K $  passes
through the point $\left( 1,0\right) $. We order the set of rays $\{
\mathcal{R}_{m,n}:| n| \leq m\leq r,m>0\} $ by increasing
slope (so that the positive slopes form the Farey sequence of order $r$),
and denote the resulting ordered sequence of rays by
$\left\{ \mathcal{L}_{j}\right\} _{j=-J}^{J}$.
Let $P_{j}$ be the intersection of $\mathcal{L}_{j}$
and $K $.

We now  define a preliminary
 domain $\cD_0$ with partially polygonal boundary, then
smooth  out the corners to get a domain $\cD$ with bounded curvature, and we
shall  then take $\Omega=\cD^*$ so that $\Omega^*=\cD^{**}=\cD$.
The boundary of the set $\cD_0$ consists of  the
polygon in the sector $\mathcal{S}=\left\{ \left( x,y\right) :0\leq \left|
y\right| \leq x\right\} $ whose edges are the segments
$\overrightarrow{P_{j}P_{j+1}}$, $-J\leq j<J$,
together with  a smooth curve of bounded  curvature
 in the closure of the complement of $\mathcal{S}$. We note that with
$P_{j}=\left( x_{j},y_{j}\right) $,  we have
\begin{equation}
1-\frac{1}{2hr^{2}}\leq x_{j}\leq 1  \label{Pyth'}
\end{equation}
by a straightforward application of  Pythagoras' theorem
and the fact that $\left| P_{j}-\left( 1,0\right)
\right| \leq 1$.

It follows that for this convex set $\cD_0$, the defining functional
$\rho_0 ^{*}$ satisfies (\ref{slim}) if $h\geq R^{1/2}
( \log R)^{-\varepsilon }$. Indeed, if $\mathcal{R}_{m,n_{\sigma }}=\mathcal{L}_{j_{\sigma }}$ for $\sigma =1,2$, then $( m,n_{\sigma }) =
\frac{m}{x_{j_{\sigma }}}\left( x_{j_{\sigma }},y_{j_{\sigma }}\right) $ and
so by (\ref{Pyth'}), we have
\begin{eqnarray*}
\left| \rho ^{*}\left( \left( m,n_{1}\right) \right) -\rho ^{*}\left( \left(
m,n_{2}\right) \right) \right|  &=&\Big|
\frac{m}{x_{j_{1}}}-\frac{m}{x_{j_{2}}}\Big| =\frac{m| x_{j_{1}}-x_{j_{2}}| }{x_{j_{1}}x_{j_{2}}} \\
&\leq &2m(hr^{2})^{-1}<2(hr)^{-1}.
\end{eqnarray*}
We now smooth out the corners to define the domain $\cD_0$.
We modify $\partial \cD_0$ in a small neighbourhood of each $P_{j}$
by inscribing a circle of radius $1$ to be
tangent to each edge incident with $P_{j}$, so that in this neighbourhood,
$\partial \cD$ is an arc $\Gamma _{j}^{*}$ of a circle of radius $1$,
where the arc $\Gamma _{j}^{*}$ is centered about the ray $\mathcal{L}_{j}$
and has diameter $\vartriangle _{j}$, where
\begin{equation} \label{farey}
c(hr^{3})^{-1}\leq \vartriangle _{j}\leq C(hr^{2})^{-1},
\end{equation}
where the proximity of $\vartriangle_j$ to the  upper or
    lower bound depends on where in the Farey sequence the slope of 
$\cL_j$ occurs. The modification of the boundary of  $\cD_0$
  that is described above is possible for
$h\geq C_{0}$, where $C_{0}$ is a
sufficiently large constant, since by the second inequality in (\ref
{farey}),
\begin{equation*}
\vartriangle _{j}
\ll r^{-2} \leq \big|
\overrightarrow{P_{j}P_{j+1}}\big| .
\end{equation*}
It is easy to see that inequality (\ref{slim}) persists for this
modification.

Now define $\Omega= \cD^*$ and let $\rho$ be the
Minkowski functional of $\Omega$.
Let $\fn^*_{(m,n)}$ be the unit normal at the boundary point of $\cD$ which lies on the ray determined by $(m,n)$.
If $\mathcal{R}_{m,n}=\mathcal{L}_{j}$
then the  arc $\Gamma_j =  (\Gamma_j^*)^*$ in $\partial \Omega$
dual to the  circular arc $\Gamma_j^*$ in $\cD$,   is itself a
 circular arc centered at
$\frac{\fn_{(m,n) }^{*}}{\rho \left( \fn_{\left( m,n\right) }^{*}\right) }$,
with curvature $1$ and diameter $\vartriangle _{j}$.
Note that the point $
\fn_{\left( m,n\right) }^{*}/
 \rho ( \fn_{(m,n) }^{*}) $ in $\partial \Omega $ has $P_{j}/|P_{j}| $
as unit normal. Thus by the first inequality in
(\ref{farey}), the cap $\mathcal{C}
(P_{j}/| P_{j}| ,\delta ) $ has diameter $\approx \delta ^{\frac{1}{2}}$ for all $0<\delta
\leq c( (hr^{3}) ^{-2}\leq c^{\prime }\vartriangle
_{j}^{2}$, and so
\begin{equation*}
(R|(m,n)|)^{1/2}\diam\big( \mathcal{C}( \tfrac{P_{m,n}}{|P_{m,n}|},
\tfrac{1}{R| ( m,n) | }) \big)
\ge c
\end{equation*}
if $(R| ( m,n) |)^{-1}\leq c(hr^{3}) ^{-2}$, and in particular if $R\geq Ch^{2}r^{6}$. 
Thus if we choose $r \le C^{1/6} (logR)^{\epsilon/3}$
and
$R^{1/2}\left( \log R\right) ^{-\varepsilon }\leq h\leq (CR)^{1/2}r^{-3},$
we have both (\ref{one}) and (\ref{slim}).\qed

\medskip

\appendix
\section{\\ Generalized distances for lattice points
in dimensions $d\ge 3$}
\label{appendix2}

Let $\rho$  be the norm associated to a   convex symmetric
 domain containing the origin, with smooth boundary and with the property that the Gaussian curvature of the boundary vanishes nowhere.
Here we are interested in lower bounds for the
distance sets
$$\Delta_{K}(E)=\{\rho \left( x-y\right) {:}x,y\in E\},$$
where $E$ will be taken as
$E(R)=\{k\in \bbZ^d, |k|\le R\}$.

This can be considered as an instance of a problem
by Erd\H os \cite{Erdos} who conjectured
for $K$ being  the unit ball for the Euclidean metric
that for {\it any } finite set $E\subset \bbR^d$ ($d\ge 2$)  one should have the estimate
\begin{equation}\label{erdos}\card (\Delta_K(E))\ge C_\eps
( \card(E))^{\frac 2d-\eps};
\end{equation}
this conjecture makes  also sense for
the more general metrics as described above and suggests the lower bound
$\card (\Delta_K(E(R)))\gc R^{2-\eps}$ in our special case,
for all metrics as described above.
The general  conjecture is open in every dimension $d\ge 2$.
For some of the best currently known  partial results and a description
of the relevant combinatorial techniques  we refer to the survey
\cite{PaSh} and other articles in the same volume.

For the case of the Euclidean metric (i.e. $K=\{|x|\le 1\}$)
the lower bound $R^{2-\eps}$ for $\Delta_K(E(R))$  is well known
and follows
from properties of the number $r(n)$ of representations of an
integer $n$  as
a sum of two squares
(see Theorems 338 and 339  in \cite{HW}). For more general metrics we
shall
 deduce the lower bound for $\card (\Delta_K(E(R)))$
from  mean discrepancy results in
\cite{ISS}
(see also the previous work by W. M\"uller \cite{M}).
Unfortunately, in two dimensions these results do not seem to
yield anything nontrivial for the distance problem.

Since the distance set
$\Delta_K(E(R))$
contains the image of $E(R)$ under $\rho$ it is
sufficient to prove the lower bound
$\card (\rho(E(R)))\ge C_\ep R^{2-\ep}$, $\ep>0$.
We have the following more precise estimate

\begin{prop}\label{rholowerbd}
Let $d\ge 3$ and let
$\Omega$ be an open convex bounded set in $\bbR^d$ containing the
origin in its interior.
Suppose that the boundary $\partial \Omega$ is  $C^\infty$, with nonvanishing Gaussian curvature.
Let $\rho$ be the Minkowski-functional associated to $\Omega$,
let
$$\bbE_R=\{k \in \bbZ^d: R/2\le  |k|\le R\}$$
and let $\rho(\bbE_R)=\{\rho(a):a\in \bbE_R\}$.
Then there exists a constant $C_0$ so that for all $R\ge C_0$ we have
$$\card(\rho(\bbE_R))\ge \begin{cases}
 R^2  &\text{ if } d\ge 4\\
R^2/\log R
&\text{ if } d=3.
\end{cases}$$
\end{prop}

\begin{proof}
Let $\alpha\ge 0$. We define for a finite set $A$ the quantity
\begin{equation*}
m_{\rho,\alpha}
(A)= \max \{\card(F): F\subset \rho(A), |s-t|>\alpha \}
\text{ for all
$s,t\in F$}\}.
\end{equation*}
In particular note that $m_{\rho,0}(A)=\rho(A)$; in fact in view
of the finiteness of $A$ we have $\rho(A)=m_{\rho,\eps_0}(A)$ for some
$\eps_0=\eps_0(A)>0$.
Moreover we let for $\eps\le 1$, $R\ge 1$
\begin{align*}S(\eps,r;A)&=\card\{k\in A; |\rho(k)-r|\le \eps\};
\\
\sigma(A,\eps)&=
\sum_{k\in A} S(2\eps,\rho(k);A).
\end{align*}
We first observe  that
for any finite set $A$  (later $A=\bbE_R$), and $\eps>0$
\begin{equation}\label{cauchyschwarz}
\card (A)\le \sqrt{m_{\rho,\eps}(A) }
\sqrt{3\sigma(A,\eps)}.
\end{equation}
Indeed if  $F$  is  a subset of $\rho(A)$ for which the
maximum in the definition of
$m_{\rho,\eps}$  is attained  then the  Cauchy-Schwarz inequality gives
\begin{equation*}
\card(A)\le \big(\card F\big)^{1/2}
\Big(\sum_{t\in F} \big( S(\eps,t; A)\big)^2 \Big)^{1/2}.
\end{equation*}
Now for fixed $k$ there are at most three intervals of the
form $[s-\eps, s+\eps]$, $s\in F$ to which $\rho(k)$ can belong.
Thus the  right hand side of the last equation is dominated by
$m_{\rho,\eps}(A)^{1/2}
(3\sum_{k\in A} S(2\eps, \rho(k);A))^{1/2}$
which yields \eqref{cauchyschwarz}.


By estimates from \cite{ISS} (namely the argument
  on p. 218/219 and the statement of Lemma 2.1 of that paper) we get for
$\eps\le R^{-1}$
\begin{align*}
\Big(R^{-d}\sum_{k\in \bbE_R} S(2\eps, \rho(k))^2
\Big)^{1/2}&\le
\Big(\sum_{k\in \bbE} S(2/R, \rho(k))^2
\Big)^{1/2}
\\
&\le C_1\Big(R^{-1}\int_{R/4}^{4R} |E(t)|^2 dt\Big)^{1/2}+ C_2 R^{d-2}
\end{align*}
where here of course $E(t)=\card(t\Omega\cap \bbZ^d)-t^d\text{vol}(\Om)$.
The results on the  mean square discrepancy in \cite{ISS} imply that
the first term is $O(R^{d-2})$ if $d\ge 4$ and  $O(R\log R)$ if $d=3$.

Now $\sigma(E_R,\eps)\le C R^d (R^{-d}\sum_{k\in \bbE_R}
S(2\eps, \rho(k))^2)^{1/2}$ and thus
 \begin{align*}
\sigma(E_R,\eps)\le C \begin{cases}
R^{2d-2} &\text{ if } d\ge 4
\\
R^{4}\log R &\text{ if } d=3
\end{cases}.
\end{align*}
Since $\card (\bbE_R)\approx R^d$  we may use \eqref{cauchyschwarz}
to obtain
for all $\eps\in (0,1/R)$ the lower bound
$m_{\rho,\eps}(A)\gc  R^2$ if $d\ge 4$ and
$m_{\rho,\eps}(A)\gc  R^2/\log R$ if $d=3$ and  this implies the
asserted lower
 bounds for the cardinality of $\rho(\bbE_R)$.
\end{proof}

{\it Remark:} For the Euclidean ball $K=B_3$
in three dimensions
one can use the mean discrepancy result by Jarn\'\i k \cite{J2}
to improve (in this very special case)
 the lower bound $R^2/\log R$ in Proposition \ref{rholowerbd}
 to $R^2/\sqrt{\log R}$.


\end{document}